\documentclass[10pt]{amsart}

\usepackage[utf8]{inputenc}         % for utf8 encoding to avoid tex accents, for example
\usepackage[T1]{fontenc}            % for T1 out encoding

\usepackage{a4wide}                 % for a4 paper size
\usepackage{amsmath, amssymb}       % various ams packages
\usepackage{amsthm}                 % ams theorem styles
\usepackage{appendix}               % appendices where you want
\usepackage{tikz}                   % graphics and diagrams
\usepackage{tikz-cd}                % commutative diagrams
\tikzcdset{ % for the style of arrow tips
arrow style=tikz
}

\usepackage{graphicx}               % check if used ??
\usepackage{stackrel}               % improved stacking of things

\usepackage{color}                  % colors for our comments, for example, and for links
\usepackage[normalem]{ulem}         % long underline across lines 
\usepackage{enumitem}               % configure appearance of lists and references
\setlist[itemize]{leftmargin=*}
\setlist[enumerate]{leftmargin=*,label=\roman*),ref=\roman*)}
\newlist{subenumerate}{enumerate}{2}
\setlist[subenumerate]{leftmargin=*,label=\alph*),ref=\alph*)}

\usepackage{blkarray}               % for arrays and tables ?? should go

\definecolor{darkblue}{rgb}{0,0,0.6} % uses color package
\usepackage[ocgcolorlinks,colorlinks=true, citecolor=darkblue, filecolor=darkblue, linkcolor=darkblue, urlcolor=darkblue]{hyperref}         % for hyperlinks
\usepackage[alphabetic]{amsrefs}

\usepackage{bbm}                    % blackboard style fonts 
\usepackage{stix}                   % just used for nsimeq 
\usepackage{eucal}                  % for Lurie's C, actually redefines the whole mathcal script
\DeclareFontFamily{T1}{cbgreek}{}
\DeclareFontShape{T1}{cbgreek}{m}{n}{<-6>  grmn0500 <6-7> grmn0600 <7-8> grmn0700 <8-9> grmn0800 <9-10> grmn0900 <10-12> grmn1000 <12-17> grmn1200 <17-> grmn1728}{}
\DeclareSymbolFont{quadratics}{T1}{cbgreek}{m}{n}
\DeclareMathSymbol{\qoppa}{\mathord}{quadratics}{19}
\DeclareMathSymbol{\Qoppa}{\mathord}{quadratics}{21}

\usepackage{comment}                % for comments that can be toggled in or out (mainly used by Fabian)
\usepackage{xspace}                 % to deal with text macros and the following trailing (or not) space. used for words additive, etc.  

\newtheoremstyle{thms}
	{}{}{\itshape}{}{\bfseries }{.}{ }
	{\thmname{#1} \thmnumber{#2}. \thmnote{\bfseries{[#3]}}}
\newtheoremstyle{thms2}
	{}{}{\itshape}{}{\bfseries }{.}{ }
	{}
\newtheoremstyle{ithreethms}
	{}{}{\itshape}{}{\bfseries }{}{ }
	{\thmname{#1} \thmnumber{#2}. \thmnote{\bfseries{[#3]}}}

\newtheoremstyle{name}
	{}{}{\itshape}{}{\bfseries }{.}{ }
	{\thmname{#1}\thmnumber{#2}\thmnote{\bfseries{[#3]}}}

\newtheoremstyle{defs}
	{}{}{\normalfont}{}{\bfseries }{.}{ }
	{\thmname{#1} \thmnumber{#2}. \thmnote{\bfseries{(#3)}}}
\newtheoremstyle{defs2}
	{}{12pt}{\normalfont}{}{\bfseries }{.}{ }
	{\thmname{#1}\thmnumber{#2}. \thmnote{\bfseries{(#3)}}}
\newtheoremstyle{ithreedefs}
	{}{}{\normalfont}{}{\bfseries }{}{ }
	{\thmname{#1} \thmnumber{#2}. \thmnote{\bfseries{(#3)}}}

\newtheoremstyle{rmk}
	{}{}{\normalfont}{}{\itshape }{.}{ }
        {}

\newtheoremstyle{claim}
	{}{}{\normalfont}{}{\itshape}{.}{ }
        {\thmname{#1} \thmnumber{#2}. \thmnote{#3}}

\theoremstyle{thms2}

\theoremstyle{thms2}

\newtheorem*{namedthm}{\namedthmname}
\newcounter{namedthm}

\makeatletter
{\endnamedthm}
\makeatother

\theoremstyle{rmk}

\theoremstyle{ithreethms}

\theoremstyle{ithreedefs}

\theoremstyle{thms2}

\swapnumbers         % so that numbers appear befor statements title, like 3.1.3 Theorem

\newcounter{rthree}
\theoremstyle{defs}
\newtheorem{definition-r-three}{Definition}[rthree]
\newtheorem{notation-r-three}[definition-r-three]{Notation}
\newtheorem{remark-r-three}[definition-r-three]{Remark}
\newtheorem{example-r-three}[definition-r-three]{Example}

\theoremstyle{thms}
\newtheorem{proposition-r-three}[definition-r-three]{Proposition}
\newtheorem{corollary-r-three}[definition-r-three]{Corollary}
\newtheorem{theorem-r-three}[definition-r-three]{Theorem}

\newcounter{rfour}
\theoremstyle{defs}
\newtheorem{definition-r-four}{Definition}[rthree]
\newtheorem{notation-r-four}[definition-r-four]{Notation}
\newtheorem{remark-r-four}[definition-r-four]{Remark}
\newtheorem{example-r-four}[definition-r-four]{Example}

\theoremstyle{thms}
\newtheorem{proposition-r-four}[definition-r-four]{Proposition}
\newtheorem{corollary-r-four}[definition-r-four]{Corollary}
\newtheorem{theorem-r-four}[definition-r-four]{Theorem}

\theoremstyle{thms}
\newtheorem{proposition}{Proposition}[section]
\newtheorem{theorem}[proposition]{Theorem}
\newtheorem{lemma}[proposition]{Lemma}
\newtheorem{corollary}[proposition]{Corollary}

\newtheorem{observation}[proposition]{Observation}

\theoremstyle{defs}
\newtheorem{definition}[proposition]{Definition}%[section]

\newtheorem{remark}[proposition]{Remark}

\theoremstyle{defs2}

\theoremstyle{rmk}

\theoremstyle{claim}

\newtheoremstyle{prf}
{}{}{\normalfont}{}{\itshape}{.}{ }
{\thmnumber{#2}\thmnote{#3}}

\theoremstyle{prf}
\newtheorem*{proofinsideproof}{Proof}

\newcommand{\spacefont}{\mathcal}                % font used for the spaces

\newcommand{\pullbacksign}{%
	\hspace{-0.25ex}%
	\tikz[baseline=(pb.base)]{%
		\draw[line width=rule_thickness, line cap=round]%
		(0,0)%
		++ (-0.5ex,-0.5ex)%
		-- ++(1ex,0ex)%
		-- ++ (0ex,1ex);%
		\filldraw (-0.3ex,0.3ex) circle (0.3pt);%
		\node (pb) at (0,0) {\phantom{x}};%
	}%
}
\newcommand{\pushoutsign}{%
	\hspace{0.25ex}%
	\tikz[baseline=(po.base)]{%
		\draw[line width=rule_thickness, line cap=round]%
		(0,0)%
		++ (-0.5ex,-0.5ex)%
		-- ++(0ex,1ex)%
		-- ++ (1ex,0ex);%
		\filldraw (0.3ex,-0.3ex) circle (0.3pt);%
		\node (po) at (0,0) {\phantom{x}};%
	}%
}
\tikzcdset{
	pullback/.style = {dash, phantom, "\pullbacksign", start anchor=center, end anchor=center},
	pushout/.style = {dash, phantom, "\pushoutsign", start anchor=center, end anchor=center}
}

\makeatletter
\newcommand{\dovl}[1]{\overline{\dbl@overline{#1}}}
\newcommand{\dbl@overline}[1]{\mathpalette\dbl@@overline{#1}}
\newcommand{\dbl@@overline}[2]{%
  \begingroup
  \sbox\z@{$\m@th#1\overline{#2}$}%
  \ht\z@=\dimexpr\ht\z@-2\dbl@adjust{#1}\relax
  \box\z@
  \ifx#1\scriptstyle\kern-\scriptspace\else
  \ifx#1\scriptscriptstyle\kern-\scriptspace\fi\fi
  \endgroup
}
\newcommand{\dbl@adjust}[1]{%
  \fontdimen8
  \ifx#1\displaystyle\textfont\else
  \ifx#1\textstyle\textfont\else
  \ifx#1\scriptstyle\scriptfont\else
  \scriptscriptfont\fi\fi\fi 3
}
\makeatother

\newcommand{\im}{\mathrm{im}}              % image
\newcommand{\Hom}{\operatorname{Hom}}            % Hom functor \warn{set Hom?}
\newcommand{\fib}{\operatorname{fib}}            % homotopy fibre
\newcommand{\cof}{\operatorname{cof}}            % homotopy quotient, cofiber

\newcommand{\colim}{\mathop{\mathrm{colim}}}     % colimits
\newcommand{\Map}{\operatorname{Map}}            % Mapping space 

\newcommand{\core}{\mathrm{Cr}}                  % groupoid core functor on $\Catx$ or $\Catp$
\newcommand{\dec}{\mathrm{dec}}                  % decalage = shift on simplicial objects

\newcommand{\nerv}{\operatorname N}              % nerve of a 1-cat
\newcommand{\Ho}{\operatorname{Ho}}              % homotopy category 

\newcommand{\op}{^\mathrm{op}}                   % opposite category
\newcommand{\Ar}{\operatorname{Ar}}              % arrow cat
\newcommand{\Twar}{\operatorname{TwAr}}          % twisted arrow cat
\newcommand{\TwR}{\operatorname{TwAr^{r}}}

\newcommand{\perf}{\mathrm{p}}                   % perfect decoration (for example in paramatrized spectra) 

\newcommand{\Mon}{\operatorname{Mon}}            % Monoid objects
\newcommand{\Grp}{\operatorname{Grp}}            % Group objects
\newcommand{\grp}{\mathrm{grp}}                  % group object decoration

\newcommand{\esd}{\operatorname{esd}}       % edgewise subdivision

\newcommand{\asscat}{\operatorname{asscat}}      % category associated to a simplicial (Segal) space, \warn{To be changed} 
\newcommand{\id}{\mathrm{id}}                % identity
\newcommand{\const}{\mathrm{const}}          % constant function

\newcommand{\Spa}{{\mathrm{Sp}}}             % spectra  \warn{To be changed}
\newcommand{\Sps}{\mathrm{An}}
\newcommand{\An}{\Sps}                % spaces  \warn{To be changed}
\newcommand{\h}{\mathrm{h}}                    % homotopy decoration
\newcommand{\Ct}{\mathrm{C_2}}                 % group $\Ct$
\newcommand{\hC}{{\h\Ct}}                      % decoration for homotopy $\Ct$ fixed points/orbits
\newcommand{\N}{\mathrm{N}}                    % norm
\newcommand{\CMonoids}{\mathrm{CMon}}            % symmetric monoidal category of comm. monoids
\newcommand{\CMon}{\CMonoids}                    % duplicate \warn{should go}

\newcommand{\Einf}{{\mathrm{E}_\infty}}          % E infinity 
\newcommand{\Eone}{{\mathrm{E}_1}}               % E one 

\newcommand{\Dperf}{{\mathcal D}^\perf}          % perfect derived category of a ring \warn{redundant}
\newcommand{\KKx}{\mathrm{KK^{\ex}}}

\newcommand{\DDelta}{{\Delta}}
\newcommand{\sAn}{\mathrm{s}\An}

\newcommand{\Cat}{\mathrm{Cat}_\infty}                  % small infinity categories
\newcommand{\Catx}{\Cat^\ex}     % small stable infinity categories
\newcommand{\Catrex}{\mathrm{Cat}_\infty^\mathrm{lex}}   % infinity categories with finite colimits and right exact functors 

\newcommand{\Catp}{\mathrm{Cat}^{\mathrm p}_\infty}     % poincare cats
\newcommand{\Span}{\operatorname{Span}}          % Span category
\newcommand{\Cob}{\mathrm{Cob}}                  % Cobordism category
\newcommand{\Seq}{\operatorname{Seq}}            % category of short exact sequences
\newcommand{\Fun}{\operatorname{Fun}}            % functor cat
\newcommand{\Nat}{\operatorname{Nat}}            % space of natural trafos
\newcommand{\nat}{\operatorname{nat}}            % spectrum of natural trafos
\newcommand{\ex}{\mathrm{st}}                    % decoration for exact functors
\newcommand{\arrowsplit}{extension splitting}

\newcommand{\Funx}{\operatorname{Fun^{ex}}}      % exact functors
\newcommand{\Null}{\operatorname{Null}}          % nullbordisms
\newcommand{\Motp}[1][]{{\if\relax\detokenize{#1}\mathrm{Mot^p}\relax\else\mathrm{Mot^p_{#1}}\fi}}  % poincare motives
\newcommand{\Motpun}[1][]{\mathrm{Mot^p_{un\if\relax\detokenize{#1}\relax\else{,}#1\fi}}}           % unstable poincare motives
\newcommand{\Motpbord}[1][]{{\if\relax\detokenize{#1}\mathrm{BMot^p}\relax\else\mathrm{BMot^p_{#1}}\fi}} % poincare-witt motives
\newcommand{\Mot}[1][]{{\if\relax\detokenize{#1}\mathrm{Mot}\relax\else\mathrm{Mot_{#1}}\fi}}            % noncommutative motives

\newcommand{\Sp}{\mathrm{Sp}}                    % symplectic group
\newcommand{\Kspace}{\operatorname{\spacefont{K}}} % K-theory space
\newcommand{\kk}{\Kspace}                        % K-theory space \warn{redundant, should go}
\newcommand{\K}{\operatorname K}                 % K-Theory \warn{spectrum or group?}
\newcommand{\KK}{\mathbb K}                      % Karoubi K-Theory (non-connective) \warn{spectrum?}

\renewcommand{\L}{\operatorname L}               % L-theory groups? \warn{Renewcommand, dangerous!}
\newcommand{\GW}{\operatorname{GW}}              % Grothendieck-Witt spectrum 
\newcommand{\Q}{\operatorname{Q}}                % Q-construction
\newcommand{\B}{\mathcal{B}}                      % cross effect \warn{wrongly used sometimes}

\newcommand{\ev}{\mathrm{ev}}                    % evaluation = double dual inclusion

\newcommand{\qshift}[1]{^{[#1]}}                 % shift (post-composition with suspension) of a hermitian category, Example: \qshift{1} 
\newcommand{\Hyp}{\operatorname{Hyp}}            % Hyperbolic categories
\newcommand{\Met}{\operatorname{Met}}            % Metabolic categories
\newcommand{\s}{\mathrm{s}}                      % same as above \warn{should go.}
\newcommand{\Poinc}{\mathrm{Pn}}                 % space of poincare objects
\newcommand{\A}{\mathcal{A}}               % an infty cat
\newcommand{\E}{\mathcal{E}}               % an infty cat (often stable)

\newcommand{\C}{\mathcal C}                % a stable infty cat
\newcommand{\D}{\mathcal{D}}               % a stable infty cat
\newcommand{\F}{\mathcal{F}}               % a functor, usually from $\Catp$
\newcommand{\QF}{\Qoppa}                   % a quadratic functor   
\usepackage{color}
\usepackage{xypic}

\author{Fabian Hebestreit}
\address{Department of Mathematics, University of Aberdeen, Scotland}
\email{f.hebestreit@abdn.ac.uk}

\author{Andrea Lachmann}
\address{Fachgruppe Mathematik und Informatik, Universit\"at Wuppertal, Germany}
\email{lachmann@uni-wuppertal.de}

\author{Wolfgang Steimle}
\address{Institut f\"ur Mathematik, Universit\"at Augsburg, Germany}
\email{wolfgang.steimle@math.uni-augsburg.de}

\title{The localisation theorem for the $K$-theory of stable $\infty$-categories}

\date{\today}

\begin{document}

\begin{abstract}
We provide a fairly self-contained account of the localisation and cofinality theorems for the algebraic $\K$-theory of stable $\infty$-categories. It is based on a general formula for the evaluation of an additive functor on a Verdier quotient closely following work of Waldhausen. We also include a new proof of the additivity theorem of $\K$-theory, strongly inspired by Ranicki's algebraic Thom construction, a short proof of the universality theorem of Blumberg, Gepner and Tabuada, and demonstrate that the cofinality theorem can be derived from the universal property alone.
\end{abstract}

\maketitle
\tableofcontents

\section{Introduction}

Algebraic $\K$-theory is one of the most prominent tools for building bridges between differential topology and number theory, ever since Quillen defined it for rings $R$. One of the simplest useful features that is immediate from his many descirptions is that $\K$-spaces of rings in fact only depend on the their categories of finitely generated projective modules and, more drastically, they are well-known to only depend on the perfect derived $\infty$-categories. This perspective was strongly advertised by Thomasson in his work on Zariski descent and makes it natural to consider algebraic $\K$-theory as a functor on the $\infty$-category of stable $\infty$-categories. That set-up also allows one to express most other important examples in a simple fashion: For example the algebraic $\K$-theory of schemes $X$ and Waldhausen's $\mathrm{A}$-theory of a spaces/animae $B$ are recovered by considering the categories of perfect complexes of quasi-coherent sheaves on $X$ and of compact spectra over $B$. Among the many existing categorical set-ups for $\K$-theory, it furthermore has the advantage of not relying on extra structure on the categories under consideration. 

In his foundational papers Quillen isolated the existence of fibre sequences relating the algebraic $\K$-spaces of rings and their localisations as a key feature of the theory, which allows one to bootstrap many calculations from the case of finite fields. While most of the fundamental theorems that go into the existence of such sequences have been transported into the setting of stable $\infty$-categories (Quillen's d\'evissage being the notable exception), there is to the best of our knowledge no account that stays within that language avoiding recourse to exact or Waldhausen categories. 
The first goal of the present note is to give such an account in particular of what has become known as the localisation theorem, i.e.\ the fact that the $\K$-space functor $\Kspace$ on stable $\infty$-categories takes Verdier sequences to fibre sequences. 
Formally, such a result first appeared in the work of Blumberg, Gepner and Tabuada on the universal property of $\mathcal K$, i.e.\ that it is the initial such functor with values in $\Einf$-groups equipped with a transformation from the core functor $\core$, which takes the maximal sub-$\infty$-groupoid of a (stable) $\infty$-category. This result sits at the heart of the trace-approach to algebraic $\K$-theory via cyclic homology and its homotopical refinements nowadays, as it facilitates easy constructions of maps out of algebraic $\K$-spaces. We give a somewhat minimalistic treatment of this fundamental result as well, in particular avoiding all mention of non-commutative motives. We also take care throughout to work with not necessarily idempotent complete categories, which allows us to provide an account of the cofinality theorem, relating the algebraic $\K$-spaces of a category to that of its idempotent completion.

All of these results rest on another basic property, called addivity, of algebraic $\K$-spaces that was isolated by Waldhausen when he extended the foundations to suit the needs of geometric topology: The additivity theorem identifies
\[\Kspace(\Ar(\C)) \simeq \Kspace(\C)^2\]
via the source and target functors $s,t \colon \Ar(\C) \to \C$; here $\Ar(\C) = \Fun([1],\C)$ is the arrow category of $\C$. 
We start our tour by providing a new direct proof of this fact for stable $\infty$-categories, that is heavily inspired by work of the first and third authors on hermitian $\K$-theory and more specifically Ranicki's algebraic Thom construction. \\

In the body of the paper an understanding of this connection (or really any aspect of hermitian K-theory or L-theory) is not required outside a few remarks but on account of the occasion, we shall take the remainder of this introduction to explain it. To this end we will for a moment assume a modicum of familiarity with Ranicki's cobordism theory for Poincar\'e chain complexes, and its reformulation in terms of Poincar\'e categories. 

Recall that in his quest to give an algebraic description of the (topological) surgery sequence, Ranicki introduced quadratic and symmetric Poincar\'e chain complexes (building on earlier work of Mishchenko) for rings with involution $R$: 
These consist of a perfect chain complex $C$ and a $\left\{\text{quadratic} \atop \text{symmetric bilinear}\right.$ form $q$ on $C$, i.e.\ an element $q$ of
\[\left\{\Omega^{\infty+d} \QF^\mathrm{q}(C) = \Hom_{\Dperf(R)}(C \otimes^\mathbb L_R C,R\qshift{-d})_\hC  \atop \Omega^{\infty+d} \QF^\s(C) = \Hom_{\Dperf(R)}(C \otimes^\mathbb L_R C,R\qshift{-d})^\hC \right.\, , \]
such that the associated polarisation $q_\sharp \colon C^{[d]} \rightarrow \mathbb R\Hom_R(C,R)$ is an equivalence. 
Here $d$ is said to be the dimension of $(C,q)$. 
Especially for $d=0$ Poincar\'e forms can be thought of as derived generalisations of $\left\{\text{quadratic} \atop \text{symmetric}\right.$ unimodular forms.
In particular Poincar\'e complexes permit a notion of Lagrangians, namely maps $f \colon L \rightarrow C$ together with a null-homotopy $\eta \colon f^*q \sim 0$, such that the induced sequence
\[L \xrightarrow{\hspace{5pt} f \hspace{5pt}} C \xrightarrow{f^* \circ q_\sharp} \mathbb R\Hom_R(L,R\qshift{-d})\]
is a fibre sequence and $\eta$ provides the null-homotopy of the composite. 
As the starting point for his theory of algebraic surgery Ranicki then proves the remarkable fact that the anima of Poincar\'e chain complexes of dimension $d$ equipped with a Lagrangian is equivalent to the anima $\left\{{\mathrm{Fm}(\Dperf(R),(\QF^\mathrm{q})\qshift{-d-1}) }\atop {\mathrm{Fm}(\Dperf(R),(\QF^\s)\qshift{-d-1})}\right.$ of chain complexes equipped with a $\left\{\text{quadratic} \atop \text{symmetric}\right.$ form of dimension $d+1$. 
On underlying objects the process takes such a $d$-dimensional Poincar\'e object $(C,q)$ with Lagrangian $L$ to the fibre of $L \rightarrow C$, and the point is that both $L$ and $C$ can be reconstructed from the induced form on the fibre. 

Now Ranicki interpreted Lagrangians in the sense above as an algebraic incarnation of null-bordisms, owing to the fact that for $M$ an oriented, compact $d$-dimensional manifold, $\mathrm{C}^*(M)$ carries a canonical $d$-dimensional symmetric Poincar\'e structure, and if $\partial W = M$ then $\mathrm C^*(W) \rightarrow \mathrm C^*(M)$ is a Lagrangian. Taking one's cue from Lefschetz duality one arrives at a notion of algebraic cobordism between two Poincar\'e objects such that the associated cobordism group is precisely Ranicki's (homotopy) $\left\{\text{quadratic} \atop \text{symmetric}\right.$ Witt- (or $\L$-) group $\left\{\L_d^\mathrm{q}(R) \atop \L_d^\s(R)\right.$. In the foundations of hermitian $\K$-theory one upgrades these cobordism groups to cobordism categories $\left\{\Cob(\Dperf(R),(\QF^\mathrm{q})\qshift{-d}) \atop \Cob(\Dperf(R),(\QF^\s)\qshift{-d})\right.$.
Then as an application of Ranicki's Thom construction one finds
\[\left\{ \Cob^\partial(\Dperf(R),(\QF^\mathrm{q})\qshift{d}) \simeq \Span(\mathrm{He}(\Dperf(R),(\QF^\mathrm{q})\qshift{d-1})) \atop \Cob^\partial(\Dperf(R),(\QF^\s)\qshift{d}) \simeq \Span(\mathrm{He}(\Dperf(R),(\QF^\s)\qshift{d-1})) \right.\, , \] 
where $\mathrm{He}$ denotes the category of forms so that $\mathrm{Fm} \simeq \core\mathrm{He}$, and by a simple cofinality argument (i.e.\ the higher categorical extension of Quillen's theorem A)
\[\left\{ |\Cob^\partial(\Dperf(R),(\QF^\mathrm{q})\qshift{d})| \atop |\Cob^\partial(\Dperf(R),(\QF^\s)\qshift{d})|\right. \simeq |\Span(\Dperf(R))| \, , \]
where $\Span(\C)$ denotes the category of spans in $\C$ (whenever $\C$ admits pullbacks).

In the present note we will follow Quillen and adopt an appropriate version of span categories as the definition of algebraic $\K$-spaces, more precisely we shall use $\Kspace(\C) = \Omega|\Span(\C)|$ as our definition (though in the stable context it is particularly easy to see that this agrees with an implementation in terms of Segal's $\mathrm S$-construction). The additivity theorem therefore takes the form
\[ \lvert \Span(\Ar(\C))\rvert \simeq {\lvert \Span(\C) \rvert} ^2.\]
Treating the left hand side as a simplistic analogue of $\Cob^\partial(\Dperf(R),\QF^\s)$ leads one looking for a similar replacement for $\Span(\mathrm{He}(\C,(\QF^\s)\qshift{-1}))$. For this end note that there is a canonical map 
\begin{align*}
	\mathrm{He}(\Dperf(R),(\QF^\s)\qshift{-1}) &\longrightarrow \TwR(\Dperf(R)) \\
	(C,q) &\longmapsto q_\sharp \, ,
\end{align*}
where $\TwR$ denotes the twisted arrow category; the superscript on the right is supposed to indicate our convention that the combined source-target map defines a right fibration $\TwR(\C) \rightarrow \C \times \C\op$. One might therefore guess that
\[\Span(\Ar(\C)) \simeq \Span(\TwR(\C))\]
for every stable $\C$ and we will show that this is indeed true by a rather simple argument. From here it is again a cofinality argument to get to the additivity theorem.

\subsection*{Notation}

From here on all categories are tacitly assumed to be $\infty$-categories. $\Catx$ denotes the category of stable categories and exact functors; $\Sp$ and $\An$ denote the categories of spectra and of $\infty$-groupoids (a.k.a.\ spaces or animae), and
\[\core\colon \Cat\to \An, \quad \Kspace\colon \Catx\to \An, \quad \K\colon \Catx\to \Sp\]
denote the groupoid-core, the space- and the spectrum-valued $\K$-theory functors.

\subsection*{Acknowledgements}

We wish to heartily thank Dustin Clausen, from whom we first learned a clean proof of the localisation theorem, for several very helpful discussions, and Christoph Winges for a very clarifying conversation about calculus of fractions. 
We further thank Ferdinand Wagner for the permission to use some of his pretty TikZ diagrams in this note.

During the preparation of this manuscript FH was a member of the Hausdorff Center for Mathematics at the University of Bonn funded by the German Research Foundation (DFG) and furthermore a member of the cluster “Mathematics Münster: Dynamics-Geometry-Structure” at the University of Münster (DFG grant nos.\ EXC 2047 390685813 and EXC 2044 390685587, resp.). FH would also like to thank the Mittag-Leffler Institute for its hospitality during the research program "Higher algebraic structures in algebra, topology and geometry", supported by the Swedish Research Council under grant no.\ 2016-06596. AL was supported by the Research Training Group "Algebro-Geometric Methods in Algebra, Arithmetic and Topology" at the University of Wuppertal (DFG grant no.\ GRK 2240) and WS by the priority program “Geometry at Infinity” (DFG grant no.\ SPP 2026) at the University of Augsburg.

\section{Localisation properties of functors on \texorpdfstring{$\Catx$}{Cat\^{}st}}
\label{sec:localisation_properties}

In the present section we briefly recall various notions of additive and localising functors, and give a brief discussion of their relation, in particular establishing a higher categorical version of Waldhausen's localisation criterion. 
We will assume the reader is familiar with Verdier sequences, that is sequences that are simultaneously fibre and cofibre sequences in $\Catx$, and refer to \cite{CDH2}*{Appendix A} for a thorough discussion.
Let us recall explicitly that a split Verdier sequence, 
i.e.\ a Verdier sequence in which both the inclusion and the projection admit both adjoints, is the same thing as a stable recollement or a semi-orthogonal decomposition into stable subcategories, see \cite{CDH2}*{Section A.2}.

For $\C$ stable, we let $\Seq(\C)$ denote the category of bifibre sequences in $\C$, i.e. the full subcategory of $\Fun([1] \times [1],\C)$ consisting of cartesian squares with lower left corner $0$. For example via
\[f \colon c \rightarrow d \quad \longmapsto \quad \begin{tikzcd} c \ar[r, "f"] \ar[d] & d \ar[d] \\ 0 \ar[r] & \cof(f)\end{tikzcd}\]
it is equivalent to $\Ar(\C) = \Fun([1],\C)$, and we shall frequently make this identification implicitly.

We very briefly recall some fundamental terminology concerning Verdier sequences:
A square
\[\begin{tikzcd} 
	\C \ar[r] \ar[d] & \C' \ar[d] \\
	\D \ar[r] & \D'
\end{tikzcd}\]
of stable categories and exact functors is a \emph{Verdier square} if it cartesian and both vertical maps are Verdier projections (i.e.\ localisations). 
It is called \emph{left} or \emph{right split} if both vertical maps are left or right split Verdier projections, respectively (i.e.\ right or left Bousfield localisations, note the order reversal). 
Finally, it is called a \emph{Karoubi square} if it becomes cartesian in the localisation of $\Catx$ at the Karoubi equivalences (i.e. those full inclusions of stable categories $\mathcal A\rightarrow \mathcal B$ which are dense in the sense that any object in $\mathcal B$ is a retract of one in $\mathcal A$), and its vertical maps are Karoubi projections, i.e.\ Verdier projections onto dense subcategories of $\D$ and $\D'$, respectively.

\begin{definition}
Let $\E$ be a category with finite limits. Then we call a functor $F \colon \Catx \rightarrow \E$ with $F(0) \simeq \, \ast$
\begin{enumerate}
\item \emph{\arrowsplit} if the combined fibre-cofibre map $(\fib,\cof) \colon \Seq(\C) \rightarrow \C^2$ induces an equivalence 
\[F(\Seq(\C)) \longrightarrow F(\C)^2\]
for every $\C \in \Catx$,
\item \emph{additive} if $F$ takes every split Verdier square to a cartesian square in $\E$, 
\item \emph{Verdier localising} if $F$ takes every Verdier square to a cartesian square in $\E$, and finally
\item \emph{Karoubi localising} if $F$ takes every Karoubi square to a cartesian square in $\E$.
\end{enumerate}
\end{definition}

From the discussion above it is hopefully obvious that the second, third and fourth condition are successively stronger. We also record:

\begin{observation}\label{observation:productspres}
All four types of functors above preserve pairwise products and moreover any group-like additive functor splits extensions.
\end{observation}

Since $\Catx$ is semi-additive it follows that any product preserving functor $F \colon \Catx \rightarrow \E$ uniquely lifts to $\Mon_\Einf(\E)$, and we call $F$ group-like if it happens to take values in $\Grp_\Einf(\E)$. 
The groupoid-core functor $\core \colon \Catx \rightarrow \An$ is an example of a non-group-like Verdier localising functor, and it particular shows that general additive functors do not split extensions.

\begin{proof}
For the preservation of products for additive functors simply note that 
\[\begin{tikzcd}
	\C \times \D \ar[r]\ar[d] & \C \ar[d] \\
	\D \ar[r] & 0
\end{tikzcd}\]
is a split Verdier square. To see that extension splitting functors preserve products, note that the map from their definition factors as
\[F(\Seq(\C)) \longrightarrow F(\C^2) \longrightarrow  F(\C)^2.\]
The composite being an equivalence implies that $F(\Seq(\C)) \longrightarrow \F(\C^2)$ admits a retraction, and it also admits a section induced by the functor $\C^2 \rightarrow \Seq(\C)$ taking $(x,y)$ to the split fibre sequence $x \rightarrow x \oplus y \rightarrow y$. 
Thus the first map in the above composition is an equivalence, and hence so is the second. But generally the map $F(\C \oplus D)  \rightarrow F(\C) \times F(\D)$ is a retract of $F((\C \oplus \D) \oplus (\C \oplus \D)) \rightarrow F(\C \oplus \D) \times F(\C \oplus \D)$, which we have just shown is an equivalence. 

The second claim follows from 
\[\begin{tikzcd}
	\C \ar[r] \ar[d] & \Seq(\C) \ar[d] \\
	0 \ar[r] & \C
\end{tikzcd}\]
being a split Verdier square together with the splitting lemma (in the category of $\Einf$-groups in $\E$); the splitting lemma itself is a direct consequence for example of \cite[Lemma 1.5.12]{CDH2}.
\end{proof}

We shall now discuss the relation between the four notions above more closely.

\subsubsection*{Additive vs.\ extension-splitting functors} 
We start with an observation from \cite{BGMN}.

\begin{definition}
Let $\KKx$ denote the ordinary category with objects stable categories and $\Hom_\KKx(\C,\D) = \K_0(\Funx(\C,\D))$ and composition induced by functor composition.
\end{definition}

Recall that $\K_0(\C)$ is defined as the quotient of the monoid $\pi_0(\core \C)$ (under direct sum) by the congruence relation generated by $x + z \sim y$ for every bifibre sequence $x \rightarrow y \rightarrow z$ in $\C$.

The natural transformation 
\[\core \Longrightarrow \pi_0\core \Longrightarrow \K_0\]
gives a functor $\Catx \rightarrow \KKx$ and we call an exact functor $\C \rightarrow \D$ a \emph{universal $\K$-equivalence} if it maps to an isomorphism in $\KKx$. 

For example if 
\[\begin{tikzcd} \C \ar[r,"f"] &  \D \ar[r,"p"] & \E \end{tikzcd}\]
is a left split Verdier sequence where $g$ and $q$ denote the left adjoints of $f$ and $p$, respectively, then the functor
\[(g,p) \colon \D \rightarrow \C \oplus \E\]
is a universal $\K$-equivalence: 
The inverse functor is $f + q \colon \C \oplus \E \rightarrow \D$, with the composition $\C \oplus \E \rightarrow \C \oplus \E$ the identity already in $\Catx$, and the other composite contained in a fibre sequence
\[qp \Longrightarrow \id_\D \Longrightarrow fg \, ,\]
see the discussion preceeding \cite[Lemma A.2.11]{CDH2}, so 
\[[(f+q) \circ (g,p)] = [fg + qp] = [\id_\D]\]
in $\K_0\Funx(\D,\D)$. The analogous claim for right split Verdier sequences holds as well.

\begin{proposition}\label{additive=splitting}
For a functor $F \colon \Catx \rightarrow \E$ with $F(0) \simeq \, \ast$, where $\E$ has finite limits, the following are equivalent:
\begin{enumerate}
\item\label{additive=splitting:item1} $F$ inverts universal $\K$-equivalences and preserves pairwise products, 
\item\label{additive=splitting:item2} $F$ is \arrowsplit, and
\item\label{additive=splitting:item3} $F$ is additive and group-like.
\end{enumerate}
\end{proposition}

The proof builds on the following classical observation of Waldhausen. We denote the functors extracting the first, second, and third entry of a fibre sequence by $\fib, m, \cof \colon \Seq(\C)\to \C$. 

\begin{lemma}\label{lem:detectionadditive}
If $F \colon \Catx \rightarrow \E$ is {\arrowsplit} then there is a canonical equivalence between $m_*$ and
\[\fib_*+\cof_* \colon F(\Seq(\C)) \longrightarrow F(\C)  \, . \] 
In particular, any {\arrowsplit} functor is group-like with the inversion map of $F(\C)$ induced by the shift functor $(-)\qshift{1} \colon \C \rightarrow \C$. 
Furthermore, a product preserving functor $F \colon \Catx \rightarrow \E$ with $F(0) \simeq 0$ is {\arrowsplit} if and only if $F$ takes
\[(s,t) \colon \Ar(\C) \longrightarrow \C^2\]
to an equivalence for every $\C \in \Catx$.
\end{lemma}

\begin{proof}
Simply note that the two functors 
\[\id \colon \Seq(\C) \longrightarrow \Seq(\C) \quad \text{and} \quad (x \rightarrow y \rightarrow z) \longmapsto (x \rightarrow x \oplus z \rightarrow z)\]
have equivalent evaluations at the first and third spot. Consequently, the latter one induces the identity on $F(\Seq(\C))$ by assumption. 
But post-composing with the evaluation at the middle term gives $\fib_* + \cof_*$ (since $\F$ preserves products by \ref{observation:productspres}).
The final item of the first claim follows from the natural bifibre sequence $x \rightarrow 0 \rightarrow x\qshift{1}$.

For the second part consider the equivalence 
\begin{align*}
	\Seq(\C) &\longrightarrow \Ar(\C) \\ 
	(x \rightarrow y \rightarrow z) &\longmapsto (x \xrightarrow{\partial} x\qshift{1})
\end{align*}
which tells us that $(\cof,\fib\qshift{1})_* \colon F\Seq(\C) \rightarrow F(\C^2) \simeq F(\C)^2$ is an equivalence. Now pre-compose with the shifting equivalence in the second factor.
\end{proof}

\begin{proof}[Proof of Proposition \ref{additive=splitting}]
The implication \ref{additive=splitting:item1} $\Rightarrow$ \ref{additive=splitting:item2} is immediate from the example preceding the statement. Conversely, being {\arrowsplit} implies that the functor
\[\mathrm hF \colon \mathrm h\Catx \rightarrow \mathrm h\E\]
factors over $\KKx$ (which immediately implies \ref{additive=splitting:item1}): 
We have to check that for a bifibre sequence $f \rightarrow g \rightarrow h$ of exact functors $\C \rightarrow \D$, say, we have $[f_*] + [h_*] = [g_*]$ in $\pi_0\Hom_\E(F(\C),F(\D))$. 
But $f,g$ and $h$ define a map $\C \rightarrow \Seq(\D)$ whence the claim follows from the previous lemma.

That \ref{additive=splitting:item3} implies \ref{additive=splitting:item2} is part of \ref{observation:productspres}, and for the final implication the previous lemma takes care of $F$ being group-like, while the example preceding the statement implies additivity.
\end{proof}

\begin{remark}
The analogous statement in the situation of Poincar\'e categories from \cites{CDH1,CDH2,CDH3}, instead of merely stable ones, is not true. 
The analogue of an {\arrowsplit} functor is a functor $F \colon \Catp \rightarrow \E$ that equates metabolic and hyperbolic Poincar\'e categories in the sense that the canonical upgrade $\Met(\C,\QF) \rightarrow \Hyp(\C)$ of the functor $(s,\cof) \colon \Ar(\C) \rightarrow \C^2$ to a Poincar\'e functor induces an equivalence for every Poincar\'e category $(\C,\QF)$. 
It is not difficult to see that such a functor inverts  universal $\GW$-equivalences (defined using the Poincar\'e refinement of $\Funx(\C,\D)$ from \cite[Section 6.2]{CDH1}).

By an argument of Schlichting a product preserving functor with this property satisfies the isotropy decomposition principle of \cite[Section 3.2]{CDH2} (i.e. \cite[Proposition 6.7 (2)]{SchlichtinghigherI} suffices as input to run the argument from \cite{CDH2}), but it need not be additive in the sense of \cite[Section 1.5]{CDH2} even if $\E$ is the category of spectra; a counterexample is the composition $\mathrm H\GW_0 \colon \Catp \rightarrow \mathcal{A}\mathrm{b} \rightarrow \Spa$. By contrast \ref{additive=splitting} or a direct check implies that $\mathrm H\K_0 \colon \Catx \rightarrow \Spa$ is additive (as are all $\mathrm H\K_i$).

The reason for this difference is that the adjoints in a split Poincar\'e-Verdier sequence are not themselves Poincar\'e functors.
\end{remark}

\subsubsection*{Additive vs.\ Verdier-localising functors}
The relation between additive and Verdier-localising functors is encapsulated by the following result, whose essence appears in Waldhausen's work as the fibration theorem \cite[Section 1.6]{waldhausen}. For a stable subcategory $\A \subseteq \B$ let $\Fun^\A(I,\B)$ denote the full subcategory of $\Fun(I,\B)$, spanned by those diagrams which take each map in $I$ to an equivalence modulo $\A$ (i.e.\ a map in $\B$ whose cofibre lies in $\A$). It is easily checked to be a stable subcategory of $\Fun(I,\B)$.

\begin{theorem}[Waldhausen]\label{Waldhausen fibration}
Given a stable subcategory $\A \subseteq \B$ the canonical maps $\const \colon \B \rightarrow \Fun^\A([n],\B)$ induce a bifibre sequence 
\[F(\A) \longrightarrow F(\B) \longrightarrow |F(\Fun^\A([-],\B))|\]
of $\Einf$-groups, whenever $F \colon \Catx \rightarrow \An$ is group-like and additive; here $F(\Fun^\A([-],\B))$ is regarded as a simplicial $\Einf$-group and the vertical bars denote its colimit.
\end{theorem}

To see the implications of this statement, recall that sifted colimits in $\Grp_\Einf(\An)$ are preserved by the forgetful functor to $\An$ so the final term is simply the geometric realisation of the simplicial anima $F(\Fun^\A([-],\B))$. Furthermore, for $M \rightarrow N \rightarrow K$ (with chosen null homotopy of the composite) being a bifibre sequence in $\Grp_\Einf(\An)$ is equivalent to the underlying sequence of anima being a fibre sequence (over the unit of $K$) and the map $\pi_0N \rightarrow \pi_0K$ being surjective.

\begin{proof}
Denote by $\dec \colon \Fun(\DDelta\op,\C) \rightarrow \Fun(\DDelta\op,\C)$ the \emph{d\'ecalage} functor induced by $[0] * - \colon \DDelta \rightarrow \DDelta$. 
Per construction the inclusions $[n] \rightarrow [0]*[n] = [1+n]$ and $[0] \rightarrow [0]*[n] = [1+n]$ induce natural transformations $\dec \Rightarrow \id$ and $\dec \Rightarrow \ev_0$. 
Recall also that $\dec(X)$ is always a split-simplicial object over $X_0$ using the latter maps and the lowest degeneracies $s_0$ of $X$ (which do not feature in the simplicial structure of $\dec(X)$); 
in particular, $X_0$ is a colimit of $\dec(X)$ by \cite[Lemma 6.1.3.16]{HTT}.

Now consider the map $d_0 \colon \Fun^\A([1+n],\B) \rightarrow \Fun^\A([n],\B)$. 
It is easily checked to be a right split Verdier projection with kernel $\A$, the requisite fully faithful right adjoint given by 
\[(b_0 \rightarrow \dots \rightarrow b_n) \longmapsto (b_0 \xrightarrow{\id} b_0 \rightarrow \dots \rightarrow b_n).\]
It follows that $\dec F \colon \Fun^\A([-],\B) \rightarrow F\Fun^\A([-],\B)$ is equifibred, i.e. that 
\[\begin{tikzcd}
F\Fun^\A([1+n],\B) \ar[r,"f^*"] \ar[d] & F\Fun^\A([1+m],\B) \ar[d] \\
F\Fun^\A([n],\B) \ar[r,"f^*"] & F\Fun^\A([m],\B)
\end{tikzcd}\]
is cartesian for every $f \colon [m] \rightarrow [n]$ in $\DDelta$: 
A square of $\Einf$-groups with right vertical map $\pi_0$-surjective (in the case at hand even split surjective) is cartesian if and only if the induced map on vertical fibres over $0$ is an equivalence. 
But this map identifies with the identity of $F(\A)$ by \ref{additive=splitting} and the analysis above. 
Note also that prior to applying $F$, the square above is \emph{not} necessarily cartesian (e.g. for $f=d_0$) , and in particular not a (right split) Verdier square.

It now follows from the equifibrancy lemma of Segal and Rezk, see e.g. \cite[Lemma 3.3.14]{CDH2} for a treatment in the present language, that 
\[\begin{tikzcd}
{|\const_{F(\A)}|}\ar[r] \ar[d] & {|\dec F\Fun^\A([-],\B)|} \ar[d] \\
{|\const_{F(0)}|} \ar[r] & {|F\Fun^\A([-],\B)|}
\end{tikzcd}\]
is cartesian, or in other words that
\[F(\A) \longrightarrow F(\B) \longrightarrow |F\Fun^\A([-],\B)|\]
is a fibre sequence. 
To finally see that it is also a cofibre sequence, note that the right hand map is (per construction) simply the inclusion of the $0$-simplices into the realisation which induces a surjection on $\pi_0$ for every simplicial anima.
\end{proof}

\begin{remark}
We will give another proof employing the relative $\Q$-construction in the final section, which passes to the Poincar\'e setting (though we will not pursue that here).
\end{remark}

Per construction, the projection $\Fun^\A([n],\B) \rightarrow \Fun([n],\B/\A)$ takes values in the subcategory spanned by those functors taking all maps in $[n]$ to equivalences in $\B/\A$. 
Since $|[n]| \simeq \, \ast$ these span the essential image of the fully faithful functor $\const \colon \B/\A \rightarrow \Fun([n],\B/\A)$, so we obtain a map 
\[|F(\Fun^\A([-],\B))| \longrightarrow F(\B/\A)\]
for every functor $F \colon \Catx \rightarrow \A$. In particular:

\begin{corollary}\label{verdiercriterion}
A grouplike additive functor $F \colon \Catx \rightarrow \An$ is Verdier localising if and only if:
\begin{enumerate}
\item the canonical map $|F(\Fun^\A([-],\B))| \longrightarrow F(\B/\A)$ constructed above is an inclusion of path components for every Verdier sequence $\A \rightarrow \B \rightarrow \B/\A$, and
\item for every Verdier square 
\[\begin{tikzcd}
\B \ar[r] \ar[d] & \B' \ar[d] \\
\B/\A \ar[r,"f"] & \B'/\A
\end{tikzcd}\]
we have
\[\mathrm{im}(\pi_0F(\B) \rightarrow \pi_0F(\B/\A)) = f^{-1}\mathrm{im}(\pi_0F(\B') \rightarrow \pi_0F(\B'/\A)) \, .\]
\end{enumerate}
\end{corollary}

For example, $|F(\Fun^\A([-],\B))| \rightarrow F(\B/\A)$ is an equivalence for every Verdier sequence if and only if $F$ is Verdier localising and $\pi_0F (\B) \rightarrow \pi_0 F(\B/\A)$ is surjective for all Verdier sequences (since $\pi_0 F\B \rightarrow \pi_0|\Fun^\A([-],\B))|$ is always surjective). We also obtain:

\begin{corollary}
An additive functor $F \colon \Catx \rightarrow \Sp$ is Verdier localising if and only if the canonical map $|F(\Fun^\A([-],\B))| \longrightarrow F(\B/\A)$ constructed above is an equivalence for every Verdier sequence $\A \rightarrow \B \rightarrow \B/\A$.
\end{corollary}

\begin{proof}
From \ref{Waldhausen fibration} it follows immediately that for additive $F$ the cofibre of $F(\A) \rightarrow F(\B)$ is given by the spectrification of $(|\Omega^\infty F(\Fun^\A([-],\B))|, |\Omega^{\infty-1}F(\Fun^\A([-],\B))|, \dots)$. But this is also a formula for the colimit of $F(\Fun^\A([-],\B))$, so we learn that 
\[F(\A) \longrightarrow F(\B) \longrightarrow |F(\Fun^\A([-],\B))|\]
is a bifibre sequence. It follows immediately that $F$ takes Verdier sequences to fibre sequences if and only $|F(\Fun^\A([-],\B))| \longrightarrow F(\B/\A)$ is an equivalence for all $\A \to \B$. But for stable targets this suffices by the following observation.
\end{proof}

\begin{lemma}\label{verdierfibresuffstable}
If $\E$ is stable, then a functor $F \colon \Catx \rightarrow \E$ with $F(0) \simeq 0$ is Verdier localising if and only if it takes Verdier sequences to (bi)fibre sequences. 
\end{lemma}

The proof is immediate from the fact that in a stable categories a commutative square is cartesian if and only if the induced map on (horizontal, say) fibres is an equivalence. 

\subsubsection*{Verdier-localising vs.\ Karoubi-localising functors}

The following is easy to check:

\begin{observation}
A functor $F \colon \Catx \rightarrow \E$ with $F(0) \simeq 0$ is Karoubi localising if and only if it is Verdier localising and inverts Karoubi equivalences.
\end{observation}

It is therefore tempting to construct a Karoubi-localising functor from a Verdier-localising one by forming
\[F \circ (-)^\natural \colon \Catx \rightarrow \E,\]
where $(-)^\natural$ denotes idempotent completion; this functor is the universal approximation of $F$ from the right by a functor inverting Karoubi equivalences. If $F$ is additive then this functor is obviously additive again, but if $F$ is Verdier-localising $F \circ (-)^\natural$ need not be so; the standard counterexample being $\mathrm{K} \colon \Catx \rightarrow \Sp$, the connective $\K$-spectrum. 

\begin{lemma}\label{idemkaroubi}
If $F \colon \Catx \rightarrow \E$ is Verdier-localising, then $F \circ (-)^\natural \colon \Catx \rightarrow \E$ is Karoubi-localising if $F$ takes pullback squares in $\Catx$, whose vertical legs are dense inclusions, to pullbacks in $\E$.
\end{lemma}

\begin{proof}
Given the previous observation it suffices to check that $F \circ (-)^\natural$ is again Verdier localising. Given then a Verdier square 
\[\begin{tikzcd} A \ar[r] \ar[d] & B \ar[d] \\
                        C \ar[r] & D\end{tikzcd}\]
with common vertical fibre $F$ consider the diagram
\[\begin{tikzcd} A^\natural \ar[r] \ar[d] & B^\natural \ar[d] \\
                        A^\natural/F \ar[r] \ar[d] & B^\natural/F \ar[d] \\
                        C^\natural \ar[r] & D^\natural.
\end{tikzcd}\]
Since idempotent completion preserves limits (e.g.\ by the comments following \ref{thomason} below) the outer square is still a pullback. Furthermore, the lower vertical maps are dense inclusions (for example since their source and target receive compatible dense inclusions from $C = A/F$ and $D=B/F$, respectively). It now follows that both squares are in fact pullbacks (in particular the upper one is a Verdier square): For the lower one it is clear that the map form $A^\natural/F$ to the pullback is fully faithful, and essential surjectivity follows from the essential surjectivity of $B^\natural \rightarrow B^\natural/F$ and it then formally follows for the upper one.

Per assumption $F$ therefore takes both individual squares to pullbacks, and consequently also the outer diagram.
\end{proof}

On stable categories there is also a dual process to idempotent completion called minimalisation, see the discussion after \ref{thomason} below. We do not know a similar criterion for functors $F \circ (-)^\mathrm{min}$, the universal approximation of $F$ from the left by a functor inverting Karoubi equivalences; such functors are again additive whenever $F$ is, but usually lose their localisation properties

\section{A reminder on \texorpdfstring{$\K$}{K}-spaces via the \texorpdfstring{$\Q$}{Q}-construction}
\label{sec:relative_Q}
The goal of the present section is to recall the $\Q$-construction, its relation to span categories and $\K$-theory. Let $\TwR(\D)$ denote the version of the twisted arrow category such that $(s,t) \colon \TwR(\D) \rightarrow \D \times \D\op$ is the right fibration classifying $\Hom_\D \colon \D\op \times \D \rightarrow \An$, see \cite[Section 5.2.1]{HA}.

\begin{definition}
For a category $\C$ with finite limits let $\Q_n(\C)$ be the full subcategory of $\Fun(\TwR([n]),\D)$ spanned by those diagrams which take every square of the form 
\begin{center}
\begin{tikzcd}[column sep=small]
	& (i \leq l)\ar[rd]\ar[ld] & \\
	(i \leq k) \ar[rd]& & (j \leq l) \ar[ld]\\
	& (j \leq k) 
\end{tikzcd}
\end{center}
to a cartesian square in $\C$.
\end{definition}

One readily checks that these categories assemble into a simplicial subcategory of $\Fun(\TwR([-]),\D)$, and in total we obtain a functor
\[\Q \colon \Catrex \rightarrow \Fun(\DDelta\op,\Catrex),\]
which we will refer to as the \emph{$\Q$-construction}. 

\begin{proposition}[Barwick]
For every category $\C$ with finite limits the simplicial category $\Q(\C)$ is a complete Segal object in $\Cat$ (and thus $\Catrex$) in the sense that the Segal maps induce an equivalences
\[\Q_n(\C) \longrightarrow \Q_1(\C) \times_{\Q_0(\C)} \dots \times_{\Q_0(\C)} \Q_1(\C)\]
and that
\[\begin{tikzcd} \Q_0(\C) \ar[r] \ar[d] & \Q_0(\C) \times \Q_0(\C) \ar[d] \\
                        \Q_3(\C) \ar[r] & \Q_1(\C) \times \Q_1(\C) \end{tikzcd}\]
is cartesian.
\end{proposition}

\begin{proof}
This precise version of the result is contained in \cite[Section 2]{HHLN2}.
\end{proof}

We can therefore apply any limit preserving functor $\Catx\to \An$ to $\Q(\C)$ and obtain a complete Segal space; recall that these are the animae which span the image of the Rezk nerve
\begin{align*}
	\nerv \colon \Catx &\longrightarrow \sAn \\ 
	\C &\longmapsto \Hom_{\Cat}([-],\C) \, .
\end{align*}
This functor $\nerv$ is fully faithful and has a left adjoint, the \emph{associated category} functor, which we will denote by
\[ \asscat \colon \sAn \longrightarrow \Catx \, . \] 

\begin{definition}
For $\C$ a category with finite limits we define 
\[\Span(\C) = \asscat(\core\Q(\C))\]
resulting in a functor $\Span \colon \Catrex \rightarrow \Cat$.
\end{definition}

We shall adopt:

\begin{definition}
For $\C$ stable we define its \emph{projective class anima} or \emph{algebraic $\K$-space} as
\[\mathcal K(\C) = \Omega|\Span(\C)|,\]
where the loop space is formed with base object $0 \in \core(\C) = \core(\Span(\C))$.
\end{definition}

In particular, the functor 
\begin{align*}
	\C &\longrightarrow \Q_1(\C) \\ 
	x &\longmapsto (0 \leftarrow x \rightarrow 0)
\end{align*} induces a map $\core(\C) \longrightarrow \kk(\C)$ natural in the input category $\C \in \Catx$. 
There results a map 
\[\pi_0(\core \C) \longrightarrow \pi_0(\kk(\C))\]
which exhibits the target as the quotient of the source by the congruence relation generated by $b \sim a + c$ whenever there is a fibre sequence $a \rightarrow b \rightarrow c$ in $\C$; this is somewhat unpleasant to see directly, but follows for example immediately from the identification of $\Q(\C)$ with the edgewise subdivision of the $\mathrm S$-construction of $\C$, compare e.g. \cite[Appendix B.1]{CDH2}. We will later use the following elementary result of Thomason, see \cite{thomason-classification}:

\begin{theorem}\label{thomason}
If $\A \rightarrow \B$ is a dense inclusion among stable categories (i.e. the functor is fully faithful and every object of $\B$ is a retract of one in $\A$), then the induced map $\K_0(\A) \rightarrow \K_0(\B)$ is injective, and sets up a bijective correspondence between dense inclusions into $\B$ (up to equivalence over $\B$) and subgroups of $\K_0(\B)$.
\end{theorem}

In particular, it follows that every stable category $\B$ admits a minimal dense stable subcategory, namely $\{b \in \B \mid [b] = 0 \in \K_0(\B)\}$. We shall denote it by $\B^\mathrm{min}$ and refer to it as the minimalisation of $\B$. In particular, $(-)^\mathrm{min} \colon \Catx \rightarrow \Catx$ is left adjoint to idempotent completion.

\section{The additivity theorem}
\label{sec:additivity_theorem}
The goal of the present section is to present a short proof of the additivity theorem for $\K$-spaces. As detailed in the previous section, we adopt $\kk(\C) \simeq \Omega|\Span(\C)|$ as the definition of these (for $\C$ stable, which is the only case we shall consider). Our goal is therefore to prove:

\begin{theorem}\label{additvitiythm}
The source and target projection give an equivalence
\[(s,t) \colon |\Span(\Ar \C)| \longrightarrow |\Span(\C)|^2\]
for every stable $\C$.
\end{theorem}

From \ref{additive=splitting} we then immediately obtain:

\begin{corollary}\label{additvitiythm2}
The functor $\kk \colon \Catx \rightarrow \An$ is additive and group-like.
\end{corollary}

\begin{proof}
Given \ref{additvitiythm} and \ref{lem:detectionadditive} note only that $\kk$ evidently preserves products, and so lifts uniquely to $\Mon_\Einf(\An)$ (as well as $\Mon_\Eone(\An)$) since $\Catx$ is semi-additive. Thus the $\Eone$-structure underlying the canonical $\Einf$-structure is also induced by the loop multiplication and this is group-like.
\end{proof}

As mentioned in the introduction the proof of \ref{additvitiythm} is strongly inspired by the algebraic Thom construction of Ranicki. Namely, we will prove the following two results which immediately imply \ref{additvitiythm}, since cofinal maps induce equivalences on realisations; see \ref{Thomconst} below for an explanation of the connection to the algebraic Thom construction.

\begin{proposition}\label{addprop1}
For $\C$ stable, there are canonical equivalences $\Span(\C) \simeq \Span(\C\op)$ and $\Span(\Ar(\C)) \simeq \Span(\TwR(\C))$ that fit together into a natural commutative square
\[\begin{tikzcd} 
\Span(\Ar(\C)) \ar[r, "\simeq"] \ar[d,"{(s,t)}"] & \Span(\TwR(\C)) \ar[d,"{(s,t)}"] \\
\Span(\C) \times \Span(\C) \ar[r,"\simeq"] & \Span(\C) \times \Span(\C\op)
\end{tikzcd}\]
\end{proposition}

It is generally true that if a category $\C$ admits pullbacks and pushouts, then $\Twar(\C)$ also has all pullbacks (so that $\Span(\TwR(\C))$ is indeed defined): Generally, the total space of a right fibration inherits pullbacks from the base by direct inspection, and this can be applied to $(s,t) \colon \TwR(\C) \rightarrow \C \times \C\op$.

\begin{proposition}\label{addprop2}
If $\C$ has finite limits and colimits and a zero object, then 
\[(s,t) \colon \Span(\TwR(\C)) \longrightarrow \Span(\C \times \C\op)\]
is cofinal.
\end{proposition}

\begin{proof}[Proof of \ref{addprop1}]
We first review the (well-known) equivalence $\Span(\C)\to\Span(\C\op)$ for stable $\C$. 
It is the identity on objects, and on morphisms it is given by first completing a span to a bicartesian square and then forgetting the initial vertex:
\[
 \begin{tikzcd}[column sep=tiny, row sep=tiny]
  & F_{01}  \ar[ld] \ar[rd] &
   & &
   & F_{01} \ar[ld]\ar[rd] & 
   &
  \\
  F_{00} && F_{11}
  & \longmapsto &
  F_{00} \ar[rd] && F_{11} \ar[ld]
  & \longmapsto &
  F_{00} \ar[rd] && F_{11} \ar[ld]
  \\
  &&
  &&
  &F_{10}&
  &&
  & F_{10}
 \end{tikzcd}
\]
where $F_{10}=F_{00}\cup_{F_{01}} F_{11}$. To define this equivalence on higher cells, denote by $\hat \Q_n(\C)\subset \Fun([n]\times [n]\op, \C)$
the full subcategory of diagrams $F$ such that each square
\[\begin{tikzcd}[column sep=tiny, row sep=tiny]
  & F_{i,j+1} \ar[ld]\ar[rd]\\
  F_{i,j} \ar[rd] && F_{i+1, j+1} \ar[ld]\\
  & F_{i+1, j}
\end{tikzcd}
\]  
is bicartesian. Restriction along the inclusion $\Twar[n]\to [n]\times [n]\op$  defines an equivalence 
\begin{equation}\label{eq:hat-Q-and-Q}\tag{$\ast$}
\hat \Q_n(\C)\longrightarrow \Q_n(\C) \, ;
\end{equation}
to see that this is indeed an equivalence consider for example the (full) subposet $\mathcal J_n \subseteq \Twar[n]$ of pairs $(i,j)$ where $i=0$ or $j=n$. Left Kan extension of functors $\mathcal J_n \rightarrow \C$ to either $\Twar[n]$ or $[n] \times [n]\op$ is then fully faithful and it follows from the pointwise formulae \cite[Proposition 4.3.2.13]{HTT} that the essential image consists exactly of $\Q_n(\C)$ or $\hat \Q_n(\C)$. Thus both these categories restrict to $\Fun(\mathcal J_n, \C)$ by an equivalence and the statement follows from 2-out-of-3 for equivalences.

The self-anti-equivalence of $[n]\times [n]\op$ that switches the entries induces an equivalence
\[\Fun([n]\times [n]\op, \C)\xrightarrow{\operatorname{op}} 
 \Fun(([n]\times [n]\op)\op, \C\op)\op \to
 \Fun([n]\times [n]\op, \C\op)\op,
\]
by pre-composition. This in turn restricts to an equivalence
\[\hat \Q_n(\C)\longrightarrow \hat \Q_n(\C\op)\op;\]
from this we obtain the desired equivalence $\Span(\C)\to \Span(\C\op)$ by applying \eqref{eq:hat-Q-and-Q} and taking groupoid cores in each simplicial degree, and then passing to the associated categories.

The equivalence $\Span(\Ar(\C))\to \Span(\Twar(\C))$ is essentially obtained by performing the above procedure in the target of the arrows; thus on morphisms it is given by the rule
\[
\begin{tikzcd}
 G_{00} \ar[d] & G_{01} \ar[l]\ar[d]\ar[r] & G_{11} \ar[d] \ar[rd, phantom, "\mapsto"] &
 G_{00} \ar[d] & G_{01} \ar[l]\ar[d]\ar[r] & G_{11} \ar[d]
 \\
 F_{00} & F_{01} \ar[l]\ar[r] & F_{11}&
 F_{00} \ar[r] & F_{10} & F_{11} \ar[l]
\end{tikzcd}
\]
On higher cells, it is obtained from a map of simplicial animae $\core\hat\Q_n(\Ar(\C))\to \core\Q_n(\Twar(\C))$ given by (restriction of) the following composite:
\begin{align*}
	\Hom_{\Cat}([n]\times [n]\op\times [1], \C) \xrightarrow{\Twar} &
	\Hom_{\Cat}(\Twar([n]\times [n]\op\times[1]), \Twar(\C)) \\
	\longrightarrow & \Hom_{\Cat}(\Twar[n], \Twar(\C))
\end{align*}
where the last map is restriction along the embedding
\begin{align*}
	\Twar[n] &\longrightarrow \Twar([n]\times [n]\op\times [1]) \\
	(i\to j) &\longmapsto ((i,j,0)\to (j,i,1)) \, .
\end{align*}
To see that this map of simplicial animae is indeed an equivalence, it suffices to consider the cases $n=0$ and $n=1$ by the Segal condition: But for $n=0$ it is the identity and for $n=1$ it is given by the rule explained above.

The commutativity of the diagram follows directly from the definition of the equivalences.
\end{proof}

\begin{remark}\label{Thomconst}
To explain the connection to Ranicki's algebraic Thom construction, let us briefly recall the latter in the language of Poincar\'e categories from \cite{CDH1}, which we shall use freely in the rest of this remark; the relation of this framework to Ranicki's categorical setting from \cite{Ranicki_bluebook} is explained in \cite[3.2.6]{HS}.

 To every Poincar\'e category $(\C,\QF)$ is associated its metabolic category $\Met(\C,\QF)$, with underlying category $\Ar(\C)$. 
Translated to this language, Ranicki showed in \cite[1.15]{Ranicki_bluebook} that associating to an arrow its fibre induces an equivalence
\[\Poinc(\Met(\C,\QF)) \rightarrow \core\mathrm{He}(\C,\QF\qshift{-1})\, ,\]
where $\Poinc$ denote the anima of Poincar\'e objects, and $\mathrm{He}$ the category of hermitian objects; see \cite[Section 2.4]{CDH1} for a treatment in the present generality. 
The name  presumably stems from the fact that for some manifold $M$ with boundary, by Atiyah duality the fibre of the metabolic object given by $\mathrm C^*(M) \rightarrow \mathrm C^*(\partial M)$ in $\C = \Dperf(\mathbb Z)$ is $C_*(\mathrm{Th} \nu_M)$, the chains of the Thom spectrum of the stable normal bundle of $M$, up to a shift.

Now for any stable category $\C$ there is another Poincar\'e category $\Hyp(\C)$ with underlying category $\C \times \C\op$. One has $\mathrm{He}(\Hyp(\C)) \simeq \TwR(\C)$ and $\Poinc(\Hyp(\C)) \simeq \core(\C)$. Using the hermitian $\Q$-construction one can then compute
\begin{align*}
\core(\Q\Ar(\C)) &\simeq \Poinc(\Hyp(\Q\Ar(\C))) \\
                        &\simeq \Poinc(\Q\Hyp\Ar(\C)) \\
                        &\simeq \Poinc(\Q\Met\Hyp(\C)) \\
                        &\simeq \Poinc(\Met(\Q\Hyp(\C))) \\
                       &\simeq \core(\mathrm{He}(\Q\Hyp(\C)\qshift{-1})) \\
                       &\simeq \core(\Q\mathrm{He}\Hyp(\C)\qshift{-1}) \\
                       &\simeq  \core(\Q\mathrm{He}\Hyp(\C)) \\
                       &\simeq \core(\Q\TwR(\C))
\end{align*}
using various commutation rules and the equivalence $\Hyp(\C)\qshift{-1} \simeq \Hyp(\C)$ via $(x,y) \mapsto (x,y\qshift{-1})$: This proves \ref{addprop1} and we found the proof given above by unwinding the effect of these equivalences. 

In \cite{CDH2} we used a similar analysis to conclude $\Poinc(\Q\Met(\C,\QF)) \simeq \core\Q\mathrm{He}(\C,\QF)$, where the left hand side defines the cobordism category $a\Cob^\partial(\C,\QF)$ of Poincar\'e objects with boundary of in $(\C,\QF)$. We then used an analogue of \ref{addprop2} to deduce $|\Cob^\partial(\C,\QF)| \simeq |\Span(\C)|$, and thus $
|\Cob(\Met(\C,\QF))| \simeq |\Cob(\Hyp(\C))|$, the hermitian analogue of \ref{additvitiythm}.
\end{remark}

For the proof of \ref{addprop2} we observe:

\begin{lemma}\label{lemmacofinal}
If $f \colon \C \rightarrow \D$ is a right fibration, $\D$ has pulbacks and $x \in \D$, then the functor $(f/x)\op \longrightarrow x/\Span(f)$ that is informally given by
\begin{align*}
	(w,f(w) \xrightarrow{\varphi} x) & \longmapsto (w,x \xleftarrow{\varphi} f(w) \xrightarrow{\id} f(w)) \\
	\intertext{admits a right adjoint given by}
	(w,x \xleftarrow{\chi} y \xrightarrow{\psi} f(w)) & \longmapsto (\bar y, f(\bar y) \simeq y \xrightarrow{\chi} x)
\end{align*}
where $\bar y \rightarrow w$ is a lift of $\psi$.
\end{lemma}

Again, we here use the observation that $\C$ inherits pullbacks from $\D$ in a fashion which is preserved by $f$ to make sense of the statement.

\begin{proof}
Unwinding definitions we have to show that
\begin{align*}
	\Hom_{x/\Span(f)}((w,x \xleftarrow{\varphi} f(w) \xrightarrow{\id} f(w)),(v,x \xleftarrow{\chi} y \xrightarrow{\psi} f(v))) \\
	\simeq \Hom_{x/f}((\bar y, f(\bar y) \simeq y \xrightarrow{\chi} x), (w,f(w) \xrightarrow{\varphi} x))
\end{align*}
naturally in $(w,f(w) \xrightarrow{\varphi} x)$. The left hand side, call it $L_\varphi$, unwinds to be the pullback
\[\begin{tikzcd}[column sep=2.3cm]
L_\varphi \ar[r] \ar[d] & \Hom_{\Span(\C)}(w,v) \ar[d] \\
\ast \ar[r,"{(x \xleftarrow{\psi}\  y \ \xrightarrow{\chi} f(v))}"] & \Hom_{\Span(\D)}(x,f(v))
\end{tikzcd}\]
where the right vertical map is 
\[\Hom_{\Span(\C)}(w,v) \xrightarrow{\Span(f)} \Hom_{\Span(\D)}(f(w),f(v)) \xrightarrow{- \circ (x \xleftarrow{\phi}\ f(w)\ \xrightarrow{id} f(w))} \Hom_{\Span(\D}(x,f(v)) \, .\]
But per definition
\[\Hom_{\Span(\C)}(w,v) \simeq \core \Fun(\TwR([1]),\C) \times_{\core(\C \times \C)} \{(w,v)\} \simeq \core(\C/w) \times_{\core \C} \core(\C/v)\]
and similarly for the lower right hand term, allowing us to rewrite this pullback as
\[\begin{tikzcd} 
L_\varphi \ar[r] \ar[d] & \core(\C/w) \times_{\core \C} \core(\C/v) \ar[d] \\
\ast \ar[r,"{(\psi,\chi)}"] & \core(\D/x) \times_{\core \D} \core(\D/f(v))
\end{tikzcd}\]
where the right vertical map identifies component-wise as
\[\C/w \xrightarrow{f} \D/f(w) \xrightarrow{ \varphi \circ -} \D/x \quad \text{and} \quad \C/v \xrightarrow{f} \D/f(v).\]
Switching the order of pullbacks, the fibre of the right hand map is contractible (since $f$ is a right fibration), so this pullback rewrites as
\[\begin{tikzcd}
L_\varphi \ar[r] \ar[d] & \core(\C/w) \times_{\core(\D/x)} \{\psi\} \ar[d] \\
\ast \ar[r,"\bar y"] & \core(\C) \times_{\core(\D)} \{y\}
\end{tikzcd}\]
where $\bar y \rightarrow v$ is the lift of $\psi$ to $\C$. Switching the order of pullbacks back this gives
\[\begin{tikzcd}
L_\varphi \ar[r] \ar[d] & \Hom_{\C}(\bar y,w) \ar[d] \\
\ast \ar[r,"\psi"] & \Hom_\D(y,x),
\end{tikzcd}\]
where the right hand vertical map is
\[\Hom_{\C}(\bar y,w) \xrightarrow{f} \Hom_{\D}(y,f(w)) \xrightarrow{\phi \circ -} \Hom_{\D}(y,x).\]
But this pullback also describes the right hand term in the equivalence we have to produce, and the whole procedure above is readily checked to be natural in $(w,\phi \colon f(w) \to x)$.
\end{proof}

\begin{proof}[Proof of \ref{addprop2}]
For the cofinality claim, we have to check that $|(c,d)/\Span(s,t)| \simeq \, \ast$ for all $(c,d) \in \C \times \C\op$. But from the lemma we find $|(c,d)/\Span(s,t)| \simeq |(s,t)/(c,d)|$ and the category $(s,t)/(c,d)$ has an initial object: One easily checks that $\id \colon 0 \rightarrow 0$ is initial in $\TwR(\C)$, and whenever a functor preserves the initial objects, all its slices inherit one.
\end{proof}

\section{The universality theorem}
\label{sec:universality_theorem}
The goal of the present section is to give a short and self-contained proof of the universal property of $\mathcal K \colon \Catx \rightarrow \An$, as first established by Blumberg, Gepner and Tabuada in \cite{BGT}. A version of the argument for higher Waldhausen categories was given by Barwick \cite{barwick2016algebraic}, and another proof in the original setting was given in \cite{CDH2}. These sources all prove a more general statement for arbitrary additive functors to $\An$ (maybe preserving filtered colimits), that gives the universal property when specified to $\core$. The following rather minimalistic argument below is adapted from \cite{steimle-note}, which proves a similar universal property in the setting of (ordinary) Waldhausen categories.

\begin{theorem}[Blumberg, Gepner, Tabuada]\label{thmbgt}
The functor $\mathcal K \colon \Catx \rightarrow \An$ is the initial group-like additive functor under $\core \colon \Catx \rightarrow \An$.
\end{theorem}

\begin{proof}
We have to show that restriction along the natural transformation $\core \Rightarrow \mathcal K$ gives an equivalence
\[\Nat(\mathcal K,F) \longrightarrow \Nat(\core,F)\]
for every grouplike additive $F \colon \Catx \rightarrow \An$. For a general functor $F  \colon \Catx \rightarrow \An$ set $G(F) = \Omega|F \Q-|$, so that $G(\core) \simeq \mathcal K$. The inclusion 
\[\C \rightarrow \Q_1(\C), \quad x \longmapsto (0 \leftarrow x \rightarrow 0)\]
induces a natural transformation $F \rightarrow G(F)$ which in turn extends to a natural transformation $\eta \colon \id \Rightarrow G$ (where $\eta_{\core}\colon \core \Rightarrow \mathcal K$ is of course the transformation considered above). Now consider the commutative diagram 
\[\begin{tikzcd} 
\Nat(G\core,F) \ar[r,"\eta_\core^*"] \ar[d,"G"] \ar[dd,bend right=80,"(\eta_F)_*","\simeq"',swap] & \Nat(\core,F) \ar[d,"G"] \ar[dd, bend left = 80, "(\eta_F)_*", "\simeq"'] \\
\Nat(GG\core,GF) \ar[r,"(G\eta_\core)^*", "\simeq"'] \ar[d,"\eta_{G\core}^*", "\simeq"'] & \Nat(G\core,GF) \ar[d,"\eta_\core^*"] \\
\Nat(G\core,GF) \ar[r,"\eta_\core^*"] & \Nat(\core,GF);
\end{tikzcd}\]
 the upper square commutes simply because $G$ is a functor, and the other three parts are consequences of the naturality of $\eta$. We now claim that the arrows labelled by $\simeq$ are equivalences: This is an immediate consequence of the following propositions and the additivity theorem (which allows us to apply the first of these to $G\core \simeq \mathcal K$). A diagram chase then implies that the entire diagram consists of equivalences.
\end{proof}

\begin{proposition}\label{etaequivgrouplike}
The transformation $\eta_F \colon F \rightarrow GF$ is an equivalence for every grouplike additive $F$.
\end{proposition}

\begin{proposition}\label{changeinauto}
The two transformations $\eta_{GF}, G\eta_F \colon GF \rightarrow GGF$ differ by an automorphism of the target. 
\end{proposition}

Proposition \ref{etaequivgrouplike} is proven in detail for example in \cite[Theorem 3.3.4]{CDH2}, but in the end the argument again goes back to Quillen and Waldhausen (a version of it is required in the proof that the iterated $\Q$- and $\mathrm{S}$- constructions define positive $\Omega$-spectra). We repeat it for completeness' sake:

\begin{definition}
For $\C$ stable we define $\Null(\C) \colon \DDelta\op \rightarrow \Catx$ as $\fib(\dec(\Q(\C)) \rightarrow \Q_0(\C))$, where the fibre is formed over $0 \in \C = \Q_0(\C)$. 
\end{definition}

Here we use the d\'ecalage functor $\dec$ from the proof of \ref{Waldhausen fibration} where $\dec(\Q(\C))\to \Q_0(\C)$ is induced by the natural transformation $\dec\Rightarrow \ev_0$. The natural transformation $\dec\rightarrow \id$ induces a natural map
\[\Null(\C) \longrightarrow \Q(\C).\]

\begin{remark}
It is not difficult to check that $\Null(\C)$ is again a complete Segal object, and that 
\[\core\Null(\C) = \nerv(0/\Span(\C)),\]
see the discussion in \cite[Section 3.3]{CDH2}.
\end{remark}

\begin{proof}[Proof of \ref{etaequivgrouplike}]
The simplicial object $\Null(\C)$ is split (over $0$) in the sense of \cite[Section 6.1.3]{HTT}, so in particular $|F\Null(\C)| \simeq \, \ast$ for every (not necessarily additive) $F \colon \Catx \rightarrow \An$.
Further it follows that the natural map
\[\eta\colon F(\C) \rightarrow \Omega|F\Q(C)|\]
is essentially by definition induced by applying $F$ and realisation to the square
\[\begin{tikzcd}
\const \C \ar[r] \ar[d] & \Null(\C) \ar[d] \\
0 \ar[r] & \Q(\C)
\end{tikzcd}\]
where the top horizontal map is induced by $\C \rightarrow \Null(\C)_0, x \mapsto (0 \leftarrow x \rightarrow 0)$. Thus we need to show that the square remains cartesian after applying $F$ and realisation. Using the equifibrancy lemma (compare again the proof of \ref{Waldhausen fibration}), we can do this by showing that the map of simplicial animae $F\Null(\C)\Rightarrow F\Q(\C)$ is equifibred for every group-like additive functor $F \colon \Catx \rightarrow \An$.

It is easy to check from the Segal condition that it suffices to treat the squares
\[\begin{tikzcd} F\Null_2(\C) \ar[r]\ar[d] & F\Q_2(\C) \ar[d] \\
                       F\Null_1(\C) \ar[r] & F\Q_1(\C)
\end{tikzcd}\]
where the vertical maps are one of $d_0,d_1$ and $d_2$. For $d_1$ and $d_2$ these squares are split Verdier squares (before applying $F$). For the remaining case we first note that $d_0$ is a split Verdier projection in both vertical maps, so that it suffices to compare vertical fibres over $0$ (since the map $\pi_0F\Q_2(\D) \rightarrow \pi_0F\Q_1(\D)$ is surjective since it is split by the degeneracy $s_0$). But on the left this fibre identifies with $\D \times \D$ and on the right with $\Ar(\D)$, and the functor between them identifies with $(d,d') \mapsto (d' \rightarrow d\oplus d')$. But this map is clearly a right inverse to $(s,\cof) \colon \Ar(\D) \rightarrow \D^2$, so an equivalence after applying $F$ by \ref{additive=splitting}. 	
\end{proof}

\begin{comment}
\begin{proof}[Proof of \ref{etaequivgrouplike}]
By the equifibrancy lemma of Segal and Rezk, see e.g. \cite[Lemma 3.3.14]{CDH2} for a treatment in the present language, it follows from the previous lemma that
\[\begin{tikzcd}
{\lvert F \const\C\rvert } \ar[r] \ar[d] & {\lvert F\Null(\C)\rvert } \ar[d] \\
{0} \ar[r] & {\lvert F\Q(\C)\rvert}
\end{tikzcd}\]
is cartesian, which is equivalent to the claim.
\end{proof}
\end{comment}

\begin{proof}[Proof of \ref{changeinauto}]
Unwinding definitions one finds that the two maps in question
\[\Omega\lvert F\Q (-)\rvert \Longrightarrow \Omega\lvert \Omega\lvert F \Q^2(-)\rvert \rvert \]
are induced by the maps into the different $\Omega$- and $\Q$-terms. In particular, the composites 
\[\Omega\lvert F\Q (-)\rvert \Longrightarrow \Omega\lvert \Omega\lvert F \Q\Q(-)\rvert \rvert \Longrightarrow \Omega^2 \lvert F\Q^2(-)\rvert \]
are exchanged by flipping both the $\Omega$ and the $\Q$-terms; here the second map is the canonical limit-colimit interchange pulling the right $\Omega$ through the outer realisation. We now claim that this is an equivalence, finishing the proof.

To this end note that for each $k \in \DDelta$ the sequence
\[ \Omega\lvert F \Q\Q_k(\C)\rvert \longrightarrow \, \ast \, \longrightarrow \lvert F \Q\Q_k(\C)\rvert \]
(with the appropriate homotopy) is not just a fibre, but also a cofibre sequence of $\Einf$-groups: 
This is equivalent to the assertion that $\pi_0\lvert F \Q\Q_k(\C)\rvert = 0$ and this holds for any stable category $\D$ in place of $\Q_k(\C)$: 
The functor 
\begin{align*}
	\Q_0\D = \D &\longrightarrow \Q_1(\D) \\ 
	x &\longmapsto (0 \leftarrow 0 \rightarrow x)
\end{align*}
composes to the identity with $d_1 \colon \Q_1(\D) \rightarrow \D$ and to $0$ with $d_0$, so the claim is a consequence of $F$ being reduced.

But then it follows that also 
\[\lvert \Omega\lvert F \Q^2(\C)\rvert \rvert \rightarrow \, \ast \, \rightarrow \lvert F \Q^2(\C)\rvert  \]
is a cofibre sequence of $\Einf$-groups, so in particular a fibre sequence of underlying animae. Looping the resulting equivalence $\lvert \Omega\lvert F \Q^2(\C)\rvert \rvert \simeq \Omega\lvert F \Q^2(\C)\rvert $ once gives the claim.
\end{proof}

\begin{remark}
The additional input needed for the more general versions of Theorem \ref{thmbgt} proved in \cite{BGT} and \cite{barwick2016algebraic} (see \cite[Section 2.7]{CDH2} for a version in the precise set-up of this paper) is that $GF = \Omega |F\Q-|$ is again additive for every additive $F \colon \Catx \rightarrow \An$. With this information the conclusion is that generally $GF$ is a group completion of $F$, i.e.\ it is initial among grouplike additive functors equipped with a transformation from $F$.

While the reader will hopefully agree that the argument we gave for the additivity of $G(\core) = \mathcal K$ in the previous section is simpler than those given in any of the references above, it does not extend beyond the case of $F = \core$, essentially because it makes use of the fact that $\core$ is defined on non-stable categories (namely $\Twar(\C)$).
\end{remark}

\section{The localisation theorem}
\label{sec:localisation_theorem}
The goal of this section is to prove:

\begin{theorem}\label{thm:K_is_localising}
The algebraic $K$-functor
\[\Kspace\colon \Catx\to \An\]
is Verdier-localising. The same is true for the spectrum-valued functor obtained from $\Kspace$ by the canonical embedding of $\Einf$-groups into spectra. 
\end{theorem}

To prove this result, we recall from Corollary \ref{verdiercriterion} that an additive and group-like functor $F\colon \Catx\to \An$ is Verdier localising provided the following condition holds:

\quad $(\ast)$ \quad \parbox{0.85\textwidth}{ %
	For any Verdier sequence $\mathcal A\to \mathcal B\to \mathcal C$, the canonical map
	$ \vert F(\Fun^{\mathcal A}([-]], \B))\vert \to F(\C) $
	is an equivalence. }

\begin{proposition}\label{prop:condition_ast_stable_under_Q}
Let $F\colon \Catx\to\An$ be an additive, space-valued functor satisfying condition $(\ast)$. Then also the functor $|F\Q|$ satisfies ($\ast$). 
\end{proposition}

In particular, it follows that if $|F\Q|$ is again additive, both functors $|F \Q|$ and $\Omega|F\Q|$ are Verdier-localising as functors $\Catx \rightarrow \Grp_\Einf(\An)$.

\begin{remark}
As mentioned previously it is generally true that $|F\Q|$ is additive whenever $F \colon \Catx \rightarrow \An$ is, but we do not prove that in this note (see \cite[Section 2.7]{CDH2}).
\end{remark}

For the proof of \ref{prop:condition_ast_stable_under_Q}, we use the following result about stable $\infty$-categories: 

\begin{lemma}
Verdier sequences are stable under applying $(-)^{\mathcal I} = \Fun(\mathcal I, -)$, for any finite poset $\mathcal I$.
\end{lemma}

\begin{proof}
Since the cotensor $(-)^{\mathcal I}$ has a left adjoint, given by the tensor $(-)_{\mathcal I}$, we deduce that $(-)^{\mathcal I}$ preserves limits (for an arbitrary category ${\mathcal I}$). On the other hand, for a finite poset ${\mathcal I}$, the functor $(-)^{\mathcal I}$ also has a right adjoint equally given by $(-)_{\mathcal I}$, see \cite[6.5.6 and 6.5.1]{CDH1}, so $(-)^{\mathcal I}$ also preserves colimits.
\end{proof}

\begin{remark}
Alternative to using the tensoring construction from \cite{CDH1}, one may prove the lemma more directly using the formulae for mapping spaces in Verdier quotients (we will recall it in \ref{thm:calculus_of_fractions} below) and in functor categories from \cite{GHNfree}, namely
\[\Nat(F,G) \simeq \lim_{f \colon x \rightarrow y \in \TwR(\mathcal J)} \Hom_\C(F(x),G(y)) \, . \]
This argument goes as follows: For $\mathcal I$ a finite poset, $\TwR(\mathcal I)$ is a finite category, so limits over it commute with filtered colimits in $\An$. It then follows easily from \ref{thm:calculus_of_fractions} that the canonical functor
\[\D^{\mathcal I}/\C^{\mathcal I} \to (\D/\C)^{\mathcal I}\]
is fully faithful whenever $\C\subseteq \D$ is a full stable subcategory. And since that formula also implies that any arrow in a Verdier quotient can be lifted with given source one can inductively show essential surjectivity: Call the length of the longest chain in $\mathcal I$ starting at some $i \in \mathcal I$ the \emph{height} $h(i)$. Given a functor $F \colon \mathcal I \rightarrow \D/\C$ and a lift $G_k \colon \mathcal I_{h \leq k} \rightarrow \D$ to the sub-poset consisting of the elements of height at most $k$, we can extend $G$ to an element $j$ of height $k+1$ by lifting $ \colim_{i < j} F(i) \rightarrow F(j)$ with source $\colim_{i<j} G(i)$. Since $I_{d \leq k+1}$ is obtained from $I_{d \leq k}$ by glueing on cones over subsets of the form $\{i \in \mathcal I\mid i<j\}$ any choices of lifts for all $j \in I$ of height $k+1$ combine into a functor $G_{k+1}$ as desired.
\end{remark}

\begin{proof}[Proof of \ref{prop:condition_ast_stable_under_Q}]
It suffices to show that the relevant map is an equivalence in each simplicial degree of the $\Q$-construction, i.e. that the canonical map
\[\vert F(\Fun^{\Q_k\A}([-], \Q_k\B))\vert \to \vert F(\Q_k\C)\vert\]
of spaces is an equivalence. By the Lemma,
\[\Q_k \A \to \Q_k \B \to \Q_k \C\]
is also a Verdier projection, in view of the equivalence
\[\Q_k\C \simeq \Fun(\mathcal J_k, \C),\]
where $\mathcal J_k\subset \Twar[k]$ is the sub-poset
\[\begin{tikzcd}[column sep=tiny, row sep=tiny]
  & 0\leq 1 \ar[ld]\ar[rd] && 1\leq 2 \ar[ld]\ar[rd]&&\ar[ld]  k-1\leq k \ar[rd] \\
0\leq 0 && 1\leq 1 && \dots && k\leq k
 \end{tikzcd}
\]
Thus, we are reduced to the case $k=0$ which holds by assumption. 
\end{proof}

To prove the first part of \ref{thm:K_is_localising}, we are left to show that the groupoid core functor $\core\colon \Catx\to \An$ satisfies condition $(\ast)$. Clearly, for any stable subcategory $\A \subseteq \B$ we have
\[\core \Fun^\A([-], \B) = \core \Fun([-], \B_\A) = \Map_{\Cat}([-], \B_\A) = \N (\B_\A)\]
where $\B_\A$ is the category of equivalences modulo $\A$ in $\B$, and $\N$ is the Rezk nerve; this category consists precisely of the maps inverted by the projection $\B \rightarrow \B/\A$ if $\A \subseteq \B$ is a Verdier inclusion, but in general this is only true up to retracts, see \cite[Lemma A.1.8]{CDH2}. Now the canonical map $\vert \N (\B_\A)\vert\to |\B_\A|$ is an equivalence, so the claim follows from:

\begin{proposition}\label{prop:category_of_equivalences_mod_A}
Let $\A\subseteq \B$ be stable subcategory and $\B_\A\subseteq \B$ the category of equivalences modulo $\A$. Then
\[\vert \B_\A\vert \to \core(\B/\A)\]
is faithful and even an equivalence if $\A\subseteq \B$ is a Verdier inclusion (i.e., closed under retracts).
\end{proposition}

\begin{remark}\label{densegivesdiscrete}
Note in particular, that the proposition implies that for dense $\C \subseteq \D$ the anima $|\D_\C|$ is discrete and in this case
\[\pi_0 \vert \D_{\C}\vert \cong \pi_0\core(\D)/\pi_0\core(\C)\cong \K_0(\D)/\K_0(\C),\]
the former by inspection, the latter by Thomason's theorem \ref{thomason}. This observation will yield the cofinality theorem in the next section.
\end{remark}

Proposition \ref{prop:category_of_equivalences_mod_A} is in turn a special case (namely, $S=\B_\A$) of the following general computation of cores in (nice enough) localisations:

\begin{proposition}\label{prop:localisation_restricts_to_localisation}
Let $S\subset \B$ be a subcategory of an $\infty$-category $\B$. Assume that the morphisms of $S$ are closed under $2$-out-of-$3$ and pushouts in $\B$. Then, the canonical functor
\[|S| = S[S^{-1}]\to \B[S^{-1}]\]
is faithful. Furthermore, the following are equivalent:
\begin{enumerate}
 \item\label{item:loc1} $|S|\subseteq \core \B[S^{-1}]$ is fully faithful.
 \item\label{item:loc2} The morphisms of $S$ satisfy  $2$-out-of-$6$ in $\B$.
 \item\label{item:loc3} A morphism of $\B$ lies in $S$ if and only if its source and target do, and it becomes invertible in $\B[S^{-1}]$.
\end{enumerate}
\end{proposition}

Here closure under pushouts means that the pushouts in $\B$ of all morphisms in $S$ exist and lie in $S$ again. The proof is based on the following formula for mapping spaces in localisations admitting a calculus of fractions, which was established by Nuiten \cite{nuiten2016localizing}. A textbook account is in {\cite[Section 7.2]{Cisinski}}. %7.2.8 \& 7.2.16

\begin{theorem}\label{thm:calculus_of_fractions}
Let $S\subset \B$ be a subcategory of an $\infty$-category $\B$. Assume that the morphisms of $S$ contain the equivalences in $\B$ and are closed under pushout in $\B$. Then, the canonical map
\[\colim_{(y\to y')\in S(y)} \Hom_{\B}(x,y') \to \Hom_{\B[S^{-1}]}(x,y)\]
is an equivalence, where $S(y)$ is the full subcategory of $y/\B$ spanned by the maps in $S$.
\end{theorem}

The special case where $\B$ is stable and $S = \B_\A$ (so that $\B[S^{-1}] = \B/\A$) was also treated by Nikolaus and Scholze in \cite[I.3.3]{NS} and enters the basic analysis of Verdier quotients. We include a short proof of the general statement for completeness' sake.

\begin{proof} 
Denote by $t_y \colon S(y) \rightarrow \B$ the functor taking targets. The proof has three steps, the first two of which in fact work for an arbitrary $S$ containing the identities of $\C$. 

\begin{enumerate}[label={\emph{Step \arabic*:}}]
	\item Whenever the functor \[\colim_{(y\to y')\in S(y)} \Hom_{\B}(-,y') \colon \B\op \rightarrow \An\] inverts $S$, it agrees with $\Hom_{\B[S^{-1}]}(-,y) \colon \B\op \rightarrow \An$.
	\item We have \[\colim_{(y\to y')\in S(y)} \Hom_{\B}(-,y') \simeq (t_y)_! \, \const_\ast \simeq \lvert -/t_y \rvert \] 
	as functors $\B\op \rightarrow \An$. 
	\item If $S$ satisfies the assumptions from the statement and $f\colon x \rightarrow x'$ is in $S$, then $f^* \colon x'/t_y \rightarrow x/t_y$ admits a left adjoint.
\end{enumerate}

\noindent Since adjunctions give equivalences on realisations, the theorem follows. 

\begin{proofinsideproof}[Proof of Step 1]
Denoting by $p \colon \B\op \rightarrow \B[S^{-1}]\op$ the localisation functor we compute
\begin{align*}p_! \left( \colim_{(y\to y')\in S(y)} \Hom_{\B}(-,y')\right) &\simeq \colim_{(y\to y')\in S(y)} p_! \Hom_{\B}(-,y') \\
	&\simeq \colim_{(y\to y')\in S(y)} \Hom_{\B[S^{-1}]}(-,y') \\
	&\simeq \Hom_{\B[S^{-1}]}(-,y) 
\end{align*}
as follows.
The first equivalence holds because left Kan extension is a left adjoint and thus preserves colimits, 
and the second one because representable functors left Kan extend to representable functors by Yoneda's lemma.
For the last equivalence observe that the functor 
\begin{align*}
	S(y) &\longrightarrow \Fun(\B\op,\An) \\ 
	(y \to y') &\longmapsto \Hom_{\B[S^{-1}]}(-,y')
\end{align*}
inverts all morphisms in $S(y)$ by two-out-of-three for equivalences in $\B[S^{-1}]$, and that $\lvert S(y) \rvert$ is contractible since it has $\id_y$ as an initial object. 

But if a functor $\B\op \rightarrow \An$ inverts $S$ (as we are assuming for the colimit), then its left Kan extension along $p$ is simply the induced functor $\B[S^{-1}]\op \rightarrow \An$ by inspection of universal properties. 
\end{proofinsideproof} 
\begin{proofinsideproof}[Proof of Step 2]
	The right hand equivalence is a direct consequence of the pointwise formula for Kan extensions:
	\[((t_y)_! \, \const_\ast)(x) \simeq \colim_{x/t_y} \const_\ast \simeq \lvert x/t_y \rvert\]
	For the left hand one we use
	\[ \colim_{f \in S(y)} \Hom_{S(y)}(g,f) \simeq \lvert g/S(y) \rvert \simeq \; \ast \] to compute
	\[(t_y)_! \, \const_* \simeq (t_y)_! \left(  \colim_{f \in S(y)} \Hom_{S(y)}(-,f) \right)\simeq \colim_{f \in S(y)} (t_y)_! \Hom_{S(y)}(-,f) \simeq \colim_{(y\to y')\in S(y)} \Hom_{\B}(-,y') \, .\]
\end{proofinsideproof}
\begin{proofinsideproof}[Proof of Step 3]
The left adjoint is easily checked to be given by taking an object $x \rightarrow z \leftarrow y$ (with left pointing arrow in $S$) to 
$x' \rightarrow p \leftarrow y$, where
\[\begin{tikzcd}
	x \ar[r, "f"] \ar[d] & x' \ar[d] \\
	z \ar[r] &p
\end{tikzcd}\]
is a pushout; this pushout exists by assumption since $f \colon x \rightarrow x'$ is in $S$ and similarly the composite $y \rightarrow z \rightarrow p$ is in $S$ since $S$ is closed under pushouts and composition. We leave it to the reader to check the adjunction property. \qedhere
\end{proofinsideproof} 
\end{proof}

\begin{remark}\label{rem:calculus_of_fractions}
Note in passing that the equivalence
\[\Hom_{\B[S^{-1}]}(x,y) \simeq \lvert x/t_y \rvert\]
arising from the proof above describes the left hand side as a certain space of zig-zags
\[x \rightarrow z \leftarrow y\]
with the left pointing arrow in $S$, as a calculus of fractions is supposed to do. We will explain in the following proof that if $S$ furthermore satisfies 2-out-of-3, then $\mathrm{Ho}(S) \subseteq \mathrm{Ho}(\B)$ indeed satisfies the classical calculus of fractions, but see \ref{calconho} below for a counterexample in general.
\end{remark}

\begin{proof}[Proof of Proposition \ref{prop:localisation_restricts_to_localisation}]
Applying Theorem \ref{thm:calculus_of_fractions} to $S\subset \B$ and $S\subset S$ (for the latter, noting that pushouts in $S$ are computed in the ambient category $\B$),  we see that
\[\Hom_{S[S^{-1}]}(x,y)\to\Hom_{\B[S^{-1}]}(x,y)\]
is computed by the formula
\[\colim_{(y\to y')\in y/S} \Hom_{S}(x,y')\to \colim_{(y\to y')\in y/S} \Hom_{\B}(x,y').\]
Since $S\subset \B$ is a subcategory, this is a directed colimit of full inclusions of subspaces, and therefore a full inclusion itself. 

This shows the first part. For the second part, we first note that all three conditions are conditions on the respective homotopy categories, and that the homotopy categories of the localisations admit a (classical) calculus of fractions if $S$ is closed under 2-out-of-3: We calculate 
\begin{align*}
\pi_0\Hom_{\B[S^{-1}]}(x,y) &= \pi_0\colim_{(y\to y')\in S(y)} \Hom_{\B}(x,y') \\
                                            &= \colim_{(y\to y')\in S(y)} \pi_0\Hom_{\B}(x,y') \\
							&= \colim_{(y\to y')\in \mathrm{Ho}(S(y))} \pi_0\Hom_{\B}(x,y') 
\end{align*}
But this colimit diagram factors through the natural functor $\mathrm{Ho}(S(y)) \rightarrow (\mathrm{Ho}(S))(y)$, which is easily checked essentially surjective and full, so it is $1$-cofinal. Thus we obtain
\[\Hom_{\mathrm{Ho}(\B[S^{-1}])}(x,y) \simeq \colim_{(y\to y')\in (\mathrm{Ho}(S))(y)} \Hom_{\mathrm{Ho}(\B)}(x,y').\] Furthermore, $(\mathrm{Ho}(S))(y)$ inherits filteredness from $S(y)$ by direct inspection and from this the axioms of a calculus of fractions \`a la Gabriel and Zisman are evident.

We now first show the implication \ref{item:loc1}$\Rightarrow$\ref{item:loc3}: If $f$ has source and target in $S$ and becomes invertible in $\B[S^{-1}]$, then under \ref{item:loc1} it is represented by a zig-zag in $\Ho(S)$ so that, by calculus of fractions and $2$-out-of-$3$, $f$ belongs itself to $S$. The implication \ref{item:loc3}$\Rightarrow$\ref{item:loc2} is trivial, since equivalences satisfy $2$-out-of-$6$. Finally, assume \ref{item:loc2} holds, and let $f$ be an invertible morphism in $\B[S^{-1}]$ between objects of $S$; we need to show that it is represented by a zigzag in $S$. 

For this, we may clearly assume that $f$ is a morphism in $\B$. 
Applying calculus of fractions again, we see that a morphism $f$ of $\B$ is split mono in the localisation if and only if, after post-composition with a morphism $g$ of $\B$, it lies in $S$. 
If $f$ is even an equivalence in the localisation, then so is $g$, and applying the same argument to $g$, we find another morphism $h$ such that $h\circ g\in S$; then $f\in S$ by $2$-out-of-$6$. 
\end{proof}

\begin{proof}[Proof of Theorem \ref{thm:K_is_localising}]
The groupoid-core functor satisfies $(\ast)$ by \ref{prop:category_of_equivalences_mod_A} and the discussion preceding it. From \ref{prop:condition_ast_stable_under_Q} we conclude that the functor $|\core \Q|$ also does, so by Corollary \ref{verdiercriterion}, $|\core \Q|\colon \Catx\to \An$ is Verdier-localising and therefore also $\Kspace=\Omega|\core\Q|$. For the second claim, we recall that a fibre sequence of $\Einf$-groups gives rise to a fibre sequence of spectra if (and only if) it is surjective on $\pi_0$, but any Verdier projection induces an epimorphism on $\K_0$ by the formula for $\K_0$ (or by noting that $\pi_0 |\core \Q(\C)|=0$ for any $\C$ and using that $|\core \Q|$ is Verdier-localising). 
\end{proof}

\begin{remark}\label{calconho}
If $S$ is closed under pushouts but not under 2-out-of-3 it need not be true that $\mathrm{Ho}(S) \subseteq \mathrm{Ho}(\B)$ satisfies the axioms for a calculus fractions as defined by Gabriel and Zisman (though something weaker is generally true, see \cite[Corollary 7.2.12]{Cisinski}). In fact, this need not even hold if $\B$ is already an ordinary category.

The following counterexample is due to Christoph Winges: Let $\B$ be an additive category and $S$ the set of inclusions of direct summands. This class is clearly closed under pushouts, but do not satisfy 2-out-of-3 unless $\B = 0$. And indeed, the two maps $0, \id_b \colon b \rightarrow b$ agree after precomposition with the map $0 \rightarrow b$ in $S$, hence they agree in $\B[S^{-1}]$ (this category in fact vanishes). But the two maps stay distinct after post-composition with any morphism in $S$ unless $b \simeq 0$. If $S(b)$ were filtered this would imply that they do not agree in $\colim_{(b\to b')\in S(a)} \Hom_{\B}(b,b')$. 

A similar example is given by taking for $\B$ the opposite of the category of Kan complexes, $S$ the trivial fibrations and any two distinct maps between contractible Kan complexes.
\end{remark}

\section{The cofinality theorem}
\label{sec:cofinality_theorem}
The first goal of this section is to formulate and prove the cofinality theorem for algebraic $K$-theory, and the second, to derive that $\Kspace$ gives rise to a Karoubi-localising functor via idempotent completion. 

The cofinality theorem follows rather directly from the methods developed for the proof of the fibration theorem. After explaining this, we give a second, independent proof of the cofinality theorem which only uses the universal property of $K$-theory and which is based on the fact that the quotient $\Einf$-monoid $\core(\A^\natural)/\core(\A)$ is group-like and discrete; recall that we denote by $(-)^\natural$ the idempotent completion of categories.  

We start by stating the cofinality theorem.  

\begin{theorem}\label{thmkcofinal}\label{thm:cofinality_for_K}
If $\A \rightarrow \B$ is a dense inclusion of stable categories, then 
\[\K_i(\A) \longrightarrow \K_i(\B)\]
is an isomorphism for $i>0$ and there is a short exact sequence
\[0 \rightarrow \K_0(\A) \rightarrow \K_0(\B) \rightarrow \pi_0\core(\B)/\pi_0\core(\A) \rightarrow 0.\] 
\end{theorem}

The statement at the level of $\K_0$ is of course part of Thomason's theorem \ref{thomason}, and we will not give an independent argument for it. One concludes from Theorem \ref{thm:cofinality_for_K} rather easily that the functor $\Kspace \colon \Catx \rightarrow \An$ satisfies the assumptions of \ref{idemkaroubi} (see e.g.\ Corollary \ref{cor:karoubian_stays_karoubian} for a generalisation), so we obtain:

\begin{corollary}\label{corkkar}
The functor $\Kspace \circ (-)^\natural \colon \Catx \rightarrow \An$ is Karoubi-localising.
\end{corollary}

By contrast, the functor $\Kspace \circ (-)^\mathrm{min}$ is not Verdier-localising: The Verdier projection $\Dperf(\mathbb Z) \rightarrow \Dperf(\mathbb Q)$ with kernel the torsion complexes, does not yield an exact sequence on $\K$-groups after minimalisation, since the map $\K_1(\Dperf(\mathbb Z)) \rightarrow \K_1(\Dperf(\mathbb Q))$ is not surjective.

Corollary \ref{corkkar} gains much of its traction from the following:

\begin{theorem}
The functor 
\[\Omega^\infty \colon \Fun(\Catx,\mathrm{Sp}) \longrightarrow \Fun(\Catx,\An)\] 
induces an equivalence between the full subcategories of Karoubi localising functors on both sides.
\end{theorem}

Versions of this result have long been known, again going back to the work of Blumberg, Gepner and Tabuada. 
The precise version above will appear as part of \cite{CDH4}, and we shall not discuss its proof any further in this note. 
It implies existence of a unique Karoubi localising functor $\KK \colon \Catx \rightarrow \mathrm{Sp}$, non-connective algebraic $K$-theory, with $\Omega^\infty \KK(\C) \simeq \Kspace(\C^\natural)$ for stable categories $\C$. 
It is this functor which is mostly used in the modern study of algebraic $K$-groups and spectra, since it (or more precisely the restriction $X \mapsto \KK(\Dperf(X))$, $\Dperf(X)$ being the perfect derived category of any scheme $X$) satisfies Zariski descent for nice enough schemes (as does any Karoubi localising functor preserving filtered colimits), while $\mathrm{K} \colon \Catx \rightarrow \mathrm{Sp}$ does not.

We now turn to the proof of Theorem \ref{thmkcofinal}.

\begin{proof}[First proof]
The first proof we give is based on the construction $\kk(\C) \simeq \Omega|\core\Q(\C)|$ and the analysis made for the localisation theorem. From \ref{Waldhausen fibration} applied to $|\core \Q-|$ we obtain a fibre sequence
\[|\core \Q\A| \longrightarrow |\core \Q\B| \longrightarrow ||\core\Q\Fun^\A([-],\B)||\]
and by inspection 
\[\core\Q\Fun^\A([-],\B) \simeq \core \Fun^{\Q\A}([-],\Q\B)\]
as bisimplicial animae. As in Section \ref{sec:localisation_theorem} we can identify 
\[|\core \Fun^{\Q_n\A}([-],\Q_n\B)| \simeq |\nerv \Q_n(\B)_{\Q_n(\A)}| \simeq |\Q_n(\B)_{\Q_n(\A)}| \, .\]
But if $\A \rightarrow \B$ is dense, so is $\Q_n(\A) \rightarrow \Q_n(\B)$. 
Thus $\Q_n(\B)/\Q_n(\A) \simeq 0$. As explained in \ref{densegivesdiscrete}, this gives that $|\Q_n(\B)_{\Q_n(\A)}|$ is discrete with components $\K_0(\Q_n(\B))/\K_0(\Q_n(\A))$. 
By direct inspection one finally finds that $\K_0(\Q(\C))$ is the edgewise subdivision of $\mathrm{Bar}(\K_0(\C))$, so that in total 
\[||\core\Q\Fun^\A([-],\B)|| \simeq |\mathrm{Bar}(\K_0(\B)/\K_0(\A))|\]
is an Eilenberg-Mac Lane space in degree $1$. Looping the fibre sequence from the start of this proof now gives the claim.
\end{proof}

The second proof of Theorem \ref{thmkcofinal} we provide rests solely on the universal property of $\mathcal K \colon \Catx \rightarrow \An$. To emphasise this, we give it in the generality of an arbitrary additive functor $F\colon \Catx\to \An$ that admits a group completion (i.e.\ an initial functor $F^\grp$ under $F$ that is group-like additive); as mentioned this is the case for any additive $F$ but we neither prove nor make use of this fact. The reader may safely consider only $F=\core$ and $F^\grp=\Kspace$ if they wish.  

\begin{definition}
We call a map $f\colon N\to M$ of $\Einf$-monoids (in $\An$) \emph{cofinal} if
\begin{enumerate}
 \item $f$ is an inclusion of a collection of path components, and 
 \item for each $x\in \pi_0(M)$ there is $x'\in \pi_0(M)$ such that $x+x'\in \pi_0(N)$.
\end{enumerate}
We call such a cofinal map \emph{dense} if in addition,
\begin{enumerate}\setcounter{enumi}{2}
 \item an element $x\in \pi_0(M)$ belongs to $\pi_0(N)$ already if there exists $y\in \pi_0(N)$ such that $x+y\in \pi_0(N)$.
\end{enumerate}
\end{definition}

The last condition is easily seen to be equivalent to the condition that the sequence of commutative monoids
\[0\to \pi_0(N)\to\pi_0(M)\to \pi_0(M)/\pi_0(N)\to 0\]
(which might generally fail to be exact in the middle) is indeed exact.

\begin{lemma}\label{lem:cofinal_submonoid}
If $f\colon N\to M$ is a cofinal map of $\Einf$-monoids, then its cofibre $M/N$ (in the category of $\Einf$-monoids) is a discrete group. 
\end{lemma}

Before proving this lemma, let us derive the cofinality theorem. Recall that every additive functor $F\colon \Catx\to \An$ automatically refines to a functor with values in $\Einf$-monoids.

\begin{definition}
We call an additive functor $F\colon \Catx\to \An$ \emph{Karoubian} if 
\begin{enumerate}
 \item every dense inclusion of stable $\infty$-categories $\A\to \B$ induces a dense map $F(\A)\to F(\B)$ of $\Einf$-monoids, and
 \item $F$ preserves pullback squares
\[\begin{tikzcd}
   \A' \ar[d] \ar[r] & \A \ar[d]\\
   \B' \ar[r] & \B
\end{tikzcd}\]
 in $\Catx$ whose (say) vertical maps are dense.
\end{enumerate}
\end{definition}

The groupoid-core functor is indeed Karoubian: The second condition holds because $\core$ commutes with all limits and for the first condition, we note that the map of $\Einf$-monoids $\core(\A)\to \core(\B)$ is clearly cofinal; furthermore if $b$ is an object in $\B$ and $a$ is an object of $\A$ such that $b\oplus a$ lies in the essential image of $\A$, then so does $b=\fib(b\oplus a\to a)$. Thus the following version of the Cofinality theorem is indeed a generalisation of \ref{thm:cofinality_for_K}.

\begin{theorem}\label{thm:cofinality}
Let $F\colon \Catx\to \An$ be an additive and Karoubian functor and $F^{\grp}$ be a group completion of $F$. For every dense inclusion $\A\to \B$ of stable $\infty$-categories, the canonical map
\[F(\B)/F(\A)\to F^{\grp}(\B)/F^\grp(\A)\]
of cofibre $\Einf$-monoids is an equivalence. Hence, $F^\grp$ induces isomorphisms
\[\pi_n F^\grp(\A)\xrightarrow\cong \pi_n F^\grp(\B), \quad n>0,\]
and a short exact sequence
\[0\to \pi_0 F^\grp(\A)\to \pi_0 F^\grp(\B) \to \pi_0 F(\B)/\pi_0 F(\A)\to 0\]
of abelian groups, where the last term denotes the quotient in the category of discrete commutative monoids.
\end{theorem}

\begin{proof}
The second statement follows from the first and Lemma \ref{lem:cofinal_submonoid}: The cofibre sequence of $\Einf$-groups
\[F^\grp(\A)\to F^\grp(\B) \to F(\B)/F(\A)\]
(with last term discrete) is a fibre sequence of animae, with last map $\pi_0$-surjective, and the functor
\[\pi_0\colon \Mon_\Einf(\An)\to \CMon\]
commutes with colimits, since it admits the discrete inclusion as a right adjoint.

Let us prove the first statement. The chain of dense inclusions
\[\A\to \B\to \A^\natural \; ( = \B^\natural)\]
induces a cofibre sequence of $\Einf$-groups
\[ F(\B)/F(\A)\to F(\A^\natural)/F(\A)\to F(\A^\natural)/F(\B) \]
and similarly with $F^{\grp}$, so it suffices to consider the case $\B=\A^\natural$. 

We claim that the functor $F'\colon \Catx\to \Mon_\Einf(\An)$ given by the formula
\[F'(\A):=F(\A^\natural)/F(\A)\]
(quotient of $\Einf$-monoids) represents the quotient $(F\circ (-)^\natural)/F$ in the category of additive functors $\Catx\to \Mon_\Einf(\An)$. To see this, it will suffice to prove that $F'$ is additive. Since $F'$ is group-like, we can combine \ref{additive=splitting} and \ref{lem:detectionadditive} with the splitting lemma to see that it suffices to prove that $F'$ sends the (split) Verdier sequence $\A \rightarrow \Ar(\A) \xrightarrow{t} \A$ to a (split) fibre sequence. In view of Lemma \ref{lem:cofinal_submonoid}, we only need to show that it induces a short exact sequences of abelian groups
\[0\to \pi_0 F'(\A)\to \pi_0 F'(\Ar(\A))\to \pi_0 F'(\A)\to 0.\]
For this, we note that short exact sequences of commutative monoids are in particular cofibre sequences of commutative monoids. Since $F$ and $F\circ (-)^\natural$ are additive, we see by commuting quotients that the sequence in question is a cofibre sequence (of commutative monoids or of abelian groups), and hence right exact. Also, the first map is (split) injective, so the sequence is indeed short exact. 

Similarly, $(F^\grp)'$ is additive and so represents the quotient of ($F^\grp\circ (-)^\natural)/F^\grp$ in the category of additive functors $\Catx\to \Grp_\Einf(\An)$: The short exact sequence
\[0\to \pi_n (F^\grp)'(\A)\to \pi_n (F^\grp)'(\B)\to \pi_n (F^\grp)'(\C)\to 0\]
is proven for $n=0$ as above, and for $n>0$ follows from the equivalence
\[\Omega (G/H)\simeq \fib(H\to G)\]
valid for any map of $\Einf$-groups. 

Next, we note that the canonical map $F\circ(-)^\natural\to F^\grp\circ(-)^\natural$ is a group completion: This follows from the fact that the endofunctor $(-)^\natural$ of $\Catx$ admits the minimalisation $(-)^\mathrm{min}$ as a left adjoint (a simple consequence of Thomason's theorem): This adjunction induces an adjunction on $\Fun(\Catx, \An)$ so we have 
equivalences
\[\nat(F\circ(-)^\natural, G) \simeq \nat(F, G\circ(-)^\mathrm{min})\simeq \nat(F^{\grp}, G\circ(-)^\mathrm{min}) \simeq \nat(F^{\grp}\circ(-)^\natural, G)\]
for every group-like additive functor $G$.

Thus, comparing universal properties, we see that the map $F'\to(F^\grp)'$ is a group completion of $F'$: But $F'$ is already group-like, so it is an equivalence.
\end{proof}

The careful reader may have noticed that the second condition of being Karoubian has not been used so far, nor has density (as opposed to cofinality) for the map induced by a dense inclusion. These extra conditions generally ensure that $F^\grp \circ (-)^\natural$ is Karoubi-localising, as we will now show. We start by observing:

\begin{corollary}\label{cor:karoubian_stays_karoubian}
If $F\colon \Catx\to \An$ is an additive and Karoubian functor, then so is any group completion $F^\grp$ of $F$.
\end{corollary}

\begin{proof}
It follows from the cofinality theorem that $F^\grp$ sends dense functors to dense maps of $\Einf$-monoids (with the last two conditions being automatic for group-like functors). It remains to prove that $F^\grp$ preserves pullback squares 
\[\begin{tikzcd}
   \A' \ar[d] \ar[r] & \A \ar[d]\\
   \B' \ar[r] & \B
\end{tikzcd}\]
 in $\Catx$ whose vertical maps are dense.
 
By assumption, this is true for $F$ and since $\pi_0 F(\A')\to \pi_0 F(\B')$ is injective, we see that also $\pi_0(F)$ sends the square to a pullback square. From density we then deduce that the map of quotient monoids
\[\pi_0 F(\B')/\pi_0 F(\A')\to \pi_0 F(\B)/\pi_0 F(\A)\]
has kernel zero. But by the cofinality theorem, this identifies with the corresponding map for $\pi_0 (F^\grp)$, so we deduce that $\pi_0(F^\grp)$ also sends the square to a pullback square. 
Applying the cofinality theorem again, we see that $F^\grp$ sends the square to a pullback square of animae. 
\end{proof}

From \ref{idemkaroubi} we now immediately obtain the following generalisation of \ref{corkkar}:

\begin{corollary}\label{cor:karoubi_localising}
Let $F\colon \Catx\to \An$ be additive and Karoubian, and let $F^\grp$ be a group completion of $F$. If $F^\grp$  is Verdier-localising, then 
\[F^\grp\circ(-)^\natural\colon \Catx\to \An\]
is Karoubi-localising.
\end{corollary}

Let us finally give the postponed

\begin{proof}[Proof of Lemma \ref{lem:cofinal_submonoid}]
We start by observing that $\pi_0(M/N) = \pi_0(M)/\pi_0(N)$ is a group by the cofinality assumption. We then need to show that 
\[\Omega(M/N)\simeq \Omega(M^\grp/N^\grp) \simeq \fib(N^\grp\to M^\grp)\]
is contractible. Since $\pi_0(N^\grp)\to \pi_0(M^\grp)$ identifies with the discrete group completion of $\pi_0(N)\to \pi_0(M)$, it is injective by cofinality, and we are left to show that the map $N^\grp\to M^\grp$ induces an equivalence on the base point components. We can prove this by showing it is a homology isomorphism because it is a map of $H$-spaces.

To show that the map on base point components is injective on homology, it suffices to prove the same for $N^\grp\to M^\grp$ itself. By the group completion theorem, this map is given by the composite
\[H_*(N)[\pi_0(N)^{-1}]\to H_*(M)[\pi_0(N)^{-1}]\to H_*(M)[\pi_0(M)^{-1}]\]
where the first map is injective by the first assumption of cofinality (since localisation is exact), and the second map is an isomorphism by the second assumption. 

To show surjectivity on the base point component, we first consider the chain of maps in the homotopy category of spaces
\[M \to M^\grp \to M^\grp_0\]
mapping an $\Einf$-monoid $M$ into the base point component of its group completion, where the first map is the canonical one and the second one subtracts $x$ in the component of $x$. On homology, these maps are given by base-change along the canonical ring homomorphisms
\[\mathbb Z[\pi_0 M] \to \mathbb Z[\pi_0 M^\grp]\to \mathbb Z\]
in view of the group completion theorem and the commutativity of $M$ for the left map and and the decomposition $M^\grp \simeq \pi_0 M^\grp \times M^\grp_0$, arising from the map $M^\grp \rightarrow M^\grp_0$ above, which is easily checked natural in the homotopy category of $\An$. We conclude the existence of a natural isomorphism
\[H_*(M^\grp_0) \cong H_*(M)/\pi_0(M)\]
where in the target we take the quotient of the monoid action in the category of graded abelian groups. 

Thus we need to show that the map $N\to M$ becomes surjective in homology after modding out the $\pi_0(M)$-action in the target. So let $a\in H_*(M)$ where we may assume that $a$ is defined in a single path component $M_x$ of $M$, for some $x\in \pi_0(M)$. By the second condition of cofinality, we may assume as well that $x\in \pi_0(N)$, in which case $a$ lifts to $H_*(N)$ by the first condition of cofinality.
\end{proof}

\section{The connection to the relative \texorpdfstring{$\Q$}{Q}-construction}
\label{sec:morerelative_Q}

Finally we briefly explain a connection between our version of Waldhausen's fibration theorem, i.e.\ Theorem \ref{Waldhausen fibration}, and the relative $\Q$-construction. In particular, this approach generalises the statement to arbitrary exact functors instead of just Verdier inclusions.

\begin{definition}
For $f \colon \C \rightarrow \D$ an exact functor between stable categories we define the \emph{relative $\Q$-construction} $\Q(f) \colon \DDelta\op \rightarrow \Catx$ by requiring
\[\begin{tikzcd}
\Q(f) \ar[d] \ar[r] & \Null(\D) \ar[d] \\
\Q(\C) \ar[r] & \Q(\D)
\end{tikzcd}\]
to be cartesian.
\end{definition}

There is a canonical map $\D \rightarrow \Q_0(f)$ whose components are given by
\begin{align*} 
	\hspace{1cm} \D &\longrightarrow \Null_0(\D) && \text{and} &\D &\longrightarrow \Q_0(\C) \hspace{2cm} \\ 
	d &\longmapsto (0 \leftarrow d \rightarrow 0) &&& d &\longmapsto 0 \, . \\
\intertext{Furthermore, the inclusion }
	&& [n] \subset [0]*[n] &\longrightarrow \TwR([0]*[n]) &&& \\ 
	&& i &\longmapsto (0_l < i_r) &&&
\end{align*}
induces a functor 
\[\Q_n(f) \longrightarrow \Fun([n],\D) \longrightarrow \Fun([n],\D/\mathrm{im}(f))\]
whose image lands in the constant functors. Since it is natural for $n \in \Delta$ we obtain a map
\[\Q(f) \Longrightarrow \const_{\D/\im(f)}.\]

We shall prove the following two statements, which together imply \ref{Waldhausen fibration} once more:

\begin{proposition}\label{Waldhausensequence}
For every group-like additive $F \colon \Catx \rightarrow \An$ and exact functor $f \colon \C \rightarrow \D$ between stable categories, the sequence
\[F(\C) \longrightarrow F(\D) \longrightarrow |F(\Q(f))|\]
is a bifibre sequence of $\Einf$-groups.
\end{proposition}

\begin{proposition}\label{SvsQ}
If $i \colon \C \rightarrow \D$ is fully faithful and exact, then there is an equivalence
\[\Fun^\C([-],\D)^{\esd} \simeq \Q(i)\]
of simplicial categories, such that
\[\begin{tikzcd} \const_\D \ar[r] \ar[d] & \ar[d] \Fun^\C([-],\D)^{\esd}\\
                        \Q(i) \ar[r]\ar[ru] & \const_{\D/\C} 
\end{tikzcd}\]
commutes.
\end{proposition}

We shall prove the two statements above in turn.

\begin{proof}[Proof of \ref{Waldhausensequence}]
By the equifibrancy lemma (see the proof of Theorem \ref{Waldhausen fibration}) the square
\[\begin{tikzcd}
{|F\Q(f)|} \ar[r] \ar[d] & {|F\Null(\D)|} \ar[d] \\
{|F\Q(\C)|} \ar[r] & {|F\Q(\D)|}
\end{tikzcd}\]
is cartesian, since the map $F \Null(\D) \to F \Q(\D)$ is equifibred as shown in the proof of \ref{etaequivgrouplike}.
In other words, 
\[|F\Q(f)| \longrightarrow |F\Q(\C)| \longrightarrow |F\Q(\D)|\]
is a fibre sequence (the case $\C=0$ gave $F(\D) \simeq \Omega|F\Q(\C)|$ in \ref{etaequivgrouplike}, since $\Q(0 \rightarrow \D) \simeq \const_\D$ by inspection). Rotating the fibre sequence above twice to the left therefore gives us the desired fibre sequence
\[F(\C) \longrightarrow F(\D) \longrightarrow |F\Q(f)| \, .\]
To see that it is also a cofibre sequence it suffices (by the long exact sequence) to check that $\pi_0|F\Q(\C)| = 0$.
But the functor
\begin{align*}
	\Q_0(\C) = \C &\longrightarrow \Q_1(\C) \\
	c &\longmapsto (0 \leftarrow 0 \rightarrow c)
\end{align*}
induces a homotopy between the $0$ and identity maps of $\pi_0|F\Q(\C)|$, which gives the claim.
\end{proof}

\begin{proof}[Proof of \ref{SvsQ}]
We shall realise both $\Q_n(i)$ and $\Fun^\C([n] * [n]\op,\D)$ as the following full subcategory $P_n$ of $\Fun(\TwR([0]*[n]*[n]\op*[0]),\D)$: A functor $F$ lies in $P_n$ if under the identification $[0]*[n]*[n]\op*[0] \cong [1+n]*[1+n]\op$,
\begin{enumerate}
\item\label{def:Pn(i)} all squares in the "left half" of $\TwR([1+n]*[1+n]\op)$ go to cartesian squares, i.e.
\begin{center}
\begin{tikzcd}[column sep=tiny]
	& F(i_l \leq m_\epsilon)\ar[rd]\ar[ld] & \\
	F(i_l \leq k_\epsilon) \ar[rd]& & F(j_l \leq m_\epsilon) \ar[ld]\\
	& F(j_l \leq k_\epsilon) 
\end{tikzcd}
\end{center}
whenever $j \leq m$ (where $\epsilon \in \{l,r\}$ and $i_\epsilon$ denotes the $i$-th element of the left respectively right factor of $[1+n]*[1+n]\op$), 
\item\label{def:Pn(ii)} $F$ vanishes on the "right half" of $\TwR([1+n]*[1+n]\op)$, i.e. we have $F(i_l \leq k_r) = 0 = F(i_r \leq k_r)$ whenever $i \geq k$, and 
\item\label{def:Pn(iii)} $F$ vanishes on the lower left corner of $\TwR([1+n]*[1+n]\op)$, i.e. $F(0_l \leq 0_l) = 0$, and 
\item\label{def:Pn(iv)} $F$ takes values in $\C$ on all spots not of the form $(0_l \leq k_\epsilon)$.
\end{enumerate}
In other words, $P_n$ consists of diagrams of the shape 
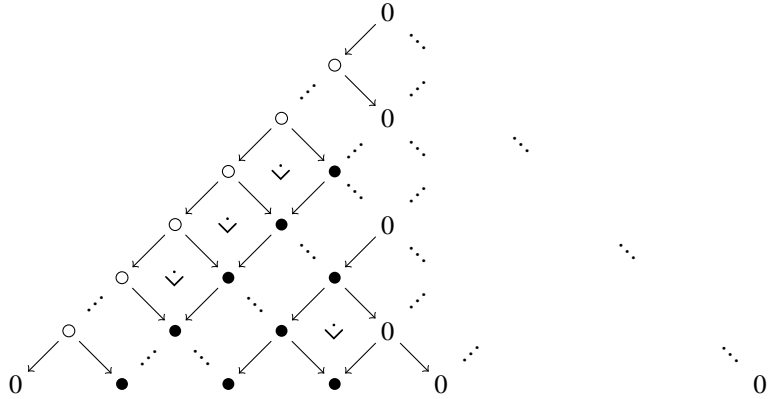
\begin{center}
	\begin{tikzpicture}[x=0.707cm,y=0.707cm,line cap=round]
		\node[inner sep=2] (00) at (0,0) {$0$};
		\draw (1,1) circle (0.5ex) coordinate (10);
		\draw (2,2) circle (0.5ex) coordinate (20);
		\draw (3,3) circle (0.5ex) coordinate (30);
		\draw (4,4) circle (0.5ex) coordinate (40);
		\draw (5,5) circle (0.5ex) coordinate (50);
		\draw (6,6) circle (0.5ex) coordinate (60);
		\node[inner sep=2] (70) at (7,7) {$0$};
		
		\fill (2,0) circle (0.5ex) coordinate (01);
		\fill (4,0) circle (0.5ex) coordinate (02);
		\fill (6,0) circle (0.5ex) coordinate (03);
		\node[inner sep=2] (04) at (8,0) {$0$};
		\node[inner sep=2] (05) at (10,0) {};
		\node[inner sep=2] (06) at (12,0) {};
		\node[inner sep=2] (07) at (14,0) {$0$};
		
		\fill (3,1) circle (0.5ex) coordinate (11);
		\fill (5,1) circle (0.5ex) coordinate (12);
		\node[inner sep=2] (13) at (7,1) {$0$};
		\node[inner sep=2] (14) at (9,1) {};
		\node[inner sep=2] (15) at (11,1) {};
		\node[inner sep=2] (16) at (13,1) {};
		\fill (4,2) circle (0.5ex) coordinate (21);
		\fill (6,2) circle (0.5ex) coordinate (22);
		\node[inner sep=2] (23) at (8,2) {};
		\node[inner sep=2] (24) at (10,2) {};
		\node[inner sep=2] (25) at (12,2) {};
		\fill (5,3) circle (0.5ex) coordinate (31);
		\node[inner sep=2] (32) at (7,3) {$0$};
		\node[inner sep=2] (33) at (9,3) {};
		\node[inner sep=2] (34) at (11,3) {};
		\fill (6,4) circle (0.5ex) coordinate (41);
		\node[inner sep=2] (42) at (8,4) {};
		\node[inner sep=2] (43) at (10,4) {};
		\node[inner sep=2] (51) at (7,5) {$0$};
		\node[inner sep=2] (52) at (9,5) {};
		\node[inner sep=2] (61) at (8,6) {};
		
		\draw[to-,shorten >=1.25ex] (00) to (10);
		\draw[to-,shorten <=1.25ex,shorten >=1.25ex] (01) to (10);
		\draw[to-,shorten <=1.25ex,shorten >=1.25ex] (02) to (12);
		\draw[to-,shorten <=1.25ex,shorten >=1.25ex] (03) to (12);
		\draw[to-,shorten <=1.25ex] (03) to (13);
		\draw[to-] (04) to (13);
		\draw[to-,shorten <=1.25ex,shorten >=1.25ex] (11) to (20);
		\draw[to-,shorten <=1.25ex,shorten >=1.25ex] (11) to (21);
		\draw[to-,shorten <=1.25ex,shorten >=1.25ex] (12) to (22);
		\draw[to-,shorten >=1.25ex] (13) to (22);
		\draw[to-,shorten <=1.25ex,shorten >=1.25ex] (20) to (30);
		\draw[to-,shorten <=1.25ex,shorten >=1.25ex] (21) to (30);
		\draw[to-,shorten <=1.25ex,shorten >=1.25ex] (21) to (31);
		\draw[to-,shorten <=1.25ex] (22) to (32);
		\draw[to-,shorten <=1.25ex,shorten >=1.25ex] (30) to (40);
		\draw[to-,shorten <=1.25ex,shorten >=1.25ex] (31) to (40);
		\draw[to-,shorten <=1.25ex,shorten >=1.25ex] (31) to (41);
		\draw[to-,shorten <=1.25ex,shorten >=1.25ex] (40) to (50);
		\draw[to-,shorten <=1.25ex,shorten >=1.25ex] (41) to (50);
		\draw[to-,shorten >=1.25ex] (51) to (60);
		\draw[to-,shorten <=1.25ex] (60) to (70);

		\path (01) to  node[pos=0.5,sloped] {$\ldots$} (11);
		\path (02) to  node[pos=0.5,sloped] {$\ldots$} (11);
		\path (04) to  node[pos=0.5,sloped] {$\ldots$} (14);
		\path (07) to  node[pos=0.5,sloped] {$\ldots$} (16);
		\path (10) to  node[pos=0.5,sloped] {$\ldots$} (20);
		\path (12) to  node[pos=0.5,sloped] {$\ldots$} (21);
		\path (13) to  node[pos=0.5,sloped] {$\ldots$} (23);
		\path (22) to  node[pos=0.5,sloped] {$\ldots$} (31);
		\path (23) to  node[pos=0.5,sloped] {$\ldots$} (32);
		\path (25) to  node[pos=0.5,sloped] {$\ldots$} (34);
		\path (32) to  node[pos=0.5,sloped] {$\ldots$} (41);
		\path (32) to  node[pos=0.5,sloped] {$\ldots$} (42);
		\path (41) to  node[pos=0.5,sloped] {$\ldots$} (51);
		\path (42) to  node[pos=0.5,sloped] {$\ldots$} (51);
		\path (43) to  node[pos=0.5,sloped] {$\ldots$} (52);
		\path (50) to  node[pos=0.5,sloped] {$\ldots$} (60);
		\path (51) to  node[pos=0.5,sloped] {$\ldots$} (61);
		\path (61) to  node[pos=0.5,sloped] {$\ldots$} (70);
		
		\node[rotate=-45] at (6,1) {$\pullbacksign$};
		\node[rotate=-45] at (3,2) {$\pullbacksign$};
		\node[rotate=-45] at (4,3) {$\pullbacksign$};
		\node[rotate=-45] at (5,4) {$\pullbacksign$};
	\end{tikzpicture}
\end{center}
where the lower left corner and the entire right half are zero, the $2n+2$ objects on the upper left diagonal marked by empty circles are in $\mathcal{D}$, the objects marked by filled circles are in the image of $\mathcal{C}$ and all squares in the left half are cartesian. 

Now consider on the one hand the inclusion
\[\alpha \colon \TwR([0]*[n]) \longrightarrow \TwR([0]*[n]*[n]\op*[0])\]
induced by the inclusion $[0] * [n] \subset ([0]*[n])*([0]*[n])\op$. Since its image lies fully in the "left half", the first and last two conditions guarantee that restriction along $\alpha$ yields a map
\[\alpha^* \colon P_n \rightarrow \Q_n(i) \, ,\]
which is clearly natural in $n \in \DDelta$.

Similarly, consider the map
\begin{align*}
	\beta \colon [n]*[n]\op &\longrightarrow \TwR([0]*[n]*[n]\op*[0])\\
	i_\epsilon &\longmapsto (0_l \leq i_\epsilon) \, .
\end{align*} The first and last conditions imply that any map in the restriction of some $F \in P_n$ goes to the pullback of some map in $\C$, so in particular to an equivalence modulo $\C$, we therefore obtain a well-defined map
\[\beta^* \colon P_n \rightarrow \Fun^\C([n]*[n]\op,\D),\]
which is again clearly natural in $n \in \Delta$. 

As mentioned the claim will now follow from both these maps $\alpha^*$ and $\beta^*$ being equivalences. We start with $\alpha^*$. An inverse is constructed by the following two-step Kan extension. First, let $T_n \subseteq \TwR([0]*[n]*[n]\op*[0])$ denote the (full) subposet given as the union of $\TwR([0]*[n])$ and the "right half" of $\TwR([0]*[n]*[n]\op*[0])$. Then consider
\[\Fun(\TwR([0]*[n]),\D) \xrightarrow{\mathrm{Lan}} \Fun(T_n,\D) \xrightarrow{\mathrm{Ran}} \Fun(\TwR([0]*[n]*[n]\op*[0]),\D).\]
As Kan extensions along fully faithful maps are fully faithful and right inverse to restriction, it remains only to check that the composite takes $\Q_n(i) \subseteq \Fun(\TwR([0]*[n]),\D)$ onto $P_n$. 
But from the pointwise formulae it is easy to see that the first Kan extension is an extension by $0$, and that the pullback condition \ref{def:Pn(i)} is equivalent to being right Kan extended from $T_n$. 
This in turn implies immediately that for $F$ satisfying conditions \ref{def:Pn(i)} and \ref{def:Pn(ii)}, condition \ref{def:Pn(iv)} is equivalent to $F\rvert_{\TwR([n])}$ taking values in $\C$, which finishes the claim (since the vanishing condition in \ref{def:Pn(iii)} is also contained in the definition of $\Q_n(i)$).

We finally treat $\beta^*$: Again the inverse is given by successive Kan extensions.  First add a zero at $(0_l \leq 0_r)$ to a functor defined on the image of $\beta$ by left Kan extension, and then right Kan extend to add zeros at $(0_l \leq 0_l)$ and the "right half", and finally left Kan extend once more to the whole of $\TwR([0]*[n]*[n]\op*[0])$. Again this process is right inverse to restriction along $\beta$ on the whole of $\Fun([n]*[n]\op,\D)$, and it remains to check that it takes $\Fun^\C([n]*[n]\op,\D)$ onto $P_n$.

But again it follows trivially from the pointwise formulae that the first two Kan extensions really are extensions by zero, precisely as required in conditions \ref{def:Pn(ii)} and \ref{def:Pn(iii)}, and that being in the image of the second left Kan extension is equivalent to the squares in condition \ref{def:Pn(i)} being pushouts, so stability implies that these are equivalent conditions. To finally see that the image of $\Fun^\C([n]*[n]\op,\D)$ also satisfies condition \ref{def:Pn(iv)}, note that (as just discussed) the value at $(i_l \leq k_\epsilon)$, with $i < k$ or $i=k$ and $\epsilon = l$, of the extension $\TwR([1+n]*[1+n]\op) \rightarrow \D$ of some $F \colon \im(\beta) \rightarrow \D$ sits in a cocartesian square
\begin{center}
\begin{tikzcd}[column sep=tiny] 
	& F(0_l \leq i_r)\ar[rd]\ar[ld] & \\
	F(0_l \leq k_\epsilon) \ar[rd]& & F(i_l \leq i_r) \ar[ld]\\
	& F(i_l \leq k_\epsilon).
\end{tikzcd}
\end{center}
But its right most term is part of the "right half" of $\TwR([0]*[n]*[n]\op*[0])$ (in fact it lies on the middle vertical line) and the upper left pointing map is an equivalence module $\C$ by assumption, so $F(i_l \leq k_\epsilon) \in \C$.
\end{proof}

\bibliographystyle{amsalpha} 
 
\renewcommand{\etalchar}[1]{$^{#1}$}
\providecommand{\bysame}{\leavevmode\hbox to3em{\hrulefill}\thinspace}
\providecommand{\MR}{\relax\ifhmode\unskip\space\fi MR }
\providecommand{\MRhref}[2]{%
  \href{http://www.ams.org/mathscinet-getitem?mr=#1}{#2}
}
\providecommand{\href}[2]{#2}

%\bibliographystyle{amsalpha}
%\ws{Check References}
%\bibliography{bib}

\end{document}